
\input amstex
\documentstyle{amsppt}
\NoBlackBoxes
\magnification=1200
\pagewidth{6.5 true  in}
\pageheight{ 8.5 true in}

\rightheadtext{}
\leftheadtext{}
\define\al{\alpha}

\define\f{\frac}
\define\n{\noindent}

\define\pa{\partial}

\topmatter
\vfill
{}\par
\title
 Deformation of spherical CR structures and the universal Picard variety
\endtitle
{}\par
\author
Jih-Hsin Cheng and I-Hsun Tsai
\endauthor
\thanks
Research supported in part by
National Science Council grants of R.O.C., NSC 86-2115M001003
and NSC 85-2121M002011, respectively.
\endthanks
\abstract
We study deformations of a spherical CR circle bundle over a Riemann
surface of genus $> 1$. Roughly speaking, there is a diffeomorphism
between such a deformation space and the unramified universal
Picard variety. On the way to parametrize the latter, we actually
give a differential-geometric proof of the structure and dimension of
the unramified universal Picard variety.
\endabstract
\endtopmatter
\document\normalbaselineskip=20pt\normalbaselines
\bigskip
\n
1. Introduction

   In this paper we study the deformation of spherical CR circle bundles
over Riemann surfaces of genus $> 1$. (for genus = 0 or 1, see [BS] for
some discussions) We find that there is a one-to-one correspondence between
such a deformation space and the so-called universal Picard variety. Let
$N$ be a closed (compact without boundary) Riemann surface of genus $g > 1$.
Let $L$ be a holomorphic line bundle over $N$ with the first Chern class
$c_1(L)$ (in $H^2(N,Z)=Z) < 0$. The universal Picard variety with
given genus $g > 1$
and $c_1 < 0$, denoted by $P_{ic}$, is the quotient space of all such
pairs $(L,N)$ modulo the equivalence relation given by holomorphic bundle
isomorphisms.

First given $(L,N)$, we can find a hermitian metric $\|\ \ \|: L\to R^+\cup
\{0\}$
such that the circle bundle $S_L\subset L$ defined by $\|\ \ \|=1$ is spherical
relative to  the induced CR structure, denoted by $J_L$ or ($H_L, J_L$).  ($H_L$
is the induced contact  bundle)(see section 2 for more details) Now fix $[(\hat L,
\hat N)]$ in $P_{ic}.$   We have the  following convention about the
regularity  of geometric objects: a  geometric object is assumed to be smooth
$(C^\infty)$  if we do not specify its regularity. We consider the deformation of
spherical CR  structures on $\hat S=S_{\hat L}$. By a theorem of Gray [Gr], we
may just fix the underlying contact bundle $\hat H=H_{\hat L}$
with the orientation induced by $\hat J=J_{\hat L}$. Let $\widetilde
{\frak S}$  denote the space of all spherical CR  manifolds $(\hat S, \hat H,
J)$ with $J$ oriented and compatible with $\hat H$. ([CL1]) Let $C_{\hat H}$
be the orientation-preserving contact diffeomorphism group
relative to $\hat H.\ \ C_{\hat H}$ acts on $\widetilde {\frak S}$ by pulling
back. Let $C^0_{\hat H}$ denote the identity component of $C_{\hat H}$. Define
the Teichmuller-type  space $\frak S^t$ to be $\widetilde
{\frak S}/C^0_{\hat H}$. Similarly we can describe $P_{ic}$ based
on a fixed background line bundle and define the Teichmuller-type space
$P_{ic}^t$. (see section 3 for  details) $P_{ic}^t$ can be
endowed with a natural complex manifold structure. (see Theorem C below)
The map $\tau:[(L,N)] \to [(S_L,H_L,J_L)]$ (equivalence relation
given by diffeomorphisms) gives rise to a map ${\tau}^t: P_{ic_0}^t
\to \frak S^t_0$. (see section 5 for definitions)

\bigskip
\n
Theorem A. (1) $\frak S^t_0$ has a  natural smooth manifold structure with
dimension equal to $8g-6.$

   (2) The map $\tau^t:  P_{ic_0}^t\to \frak S^t_0$ is a
diffeomorphism.

\bigskip
   Theorem A is in the same spirit as that of describing Teichmuller space
by conformal classes. It is known in Teichmuller theory that we can pick
up a unique hyperbolic metric as a representative for each conformal class. The
similar situation occurs for our  spherical CR manifolds. Let $\frak M_{-1,0}$
denote the quotient space of all  pseudohermitian manifolds $(M,H,J, \theta)$
with $(M,H)$ being contact-diffeomorphic to $(\hat S, \hat H)$ so that the
(pseudohermitian or Tanaka-Webster) curvature $R_{J,\theta}$ equals -1 and the torsion
$A_{J,\theta}$ vanishes modulo the equivalence relation given by
diffeomorphisms. ([We1],[Tan]) It follows that such $(M,H,J)$ is spherical and
for $(\hat S, \hat H, \hat J=J_{\hat L})$ we can always pick up a unique
contact form $\hat \theta = \theta_{\hat L}$ with $R_{\hat J,\hat\theta}=-1$
and $A_{\hat J,\hat\theta}=0.\ \ \hat H$ is given a natural orientation by
claiming $(v, \hat Jv$) is an oriented basis of $\hat H$ for any nonzero $v$ in
$\hat H$. A pseudohermitian structure $(J,\theta)$ on $(\hat S, \hat H$) is
called oriented if both J and $\theta$ are oriented for $\hat H$ . ([CL1])
To study $\frak M_{-1,0}$ we may just fix $(M, H)=(\hat S, \hat H)$ and
consider  the space of all oriented pseudohermitian structures $(J, \theta)$ on
$(\hat S, \hat H)$ with $R_{J,\theta}=-1, A_{J,\theta}=0$, denoted by $\widetilde
{\frak M}_{-1,0}$.  It is clear that $C_{\hat H}$ acts on $\widetilde{\frak
M}_{-1,0}$ by pulling back and ${\frak M}_{-1,0}= \widetilde {\frak
M}_{-1,0}/C_{\hat H}$.  Endow $\widetilde {\frak M}_{-1,0}$ with the $C^\infty$
topology and ${\frak M}_{-1,0}$ with the quotient topology.  Let $\widetilde
{\frak M}^0_{-1,0}$ be the connected component of $\widetilde {\frak
M}_{-1,0},$ containing $(\hat J,\hat \theta)$. Define the Teichmuller-type
space ${\frak T}_{-1,0}$ to be $\widetilde {\frak M}^0_{-1,0}/C^0_{\hat H}$.

\bigskip
Corollary B. The map $\iota: \frak T_{-1,0}\to \frak S^t_0$ given by
$\iota[(\hat S, \hat H, J, \theta)] = [(\hat S, \hat H, J)]$ is well
defined and a homeomorphism.

\bigskip
  Thus we can endow $\frak T_{-1,0}$ with the smooth manifold structure induced
from $\frak S^t_0$ through $\iota$.

   The universal Picard variety (or Jacobian variety) plays an important
role for many problems in algebraic geometry. Thus our differential-geometric
proof of the structure and dimension of the unramified universal Picard
variety $P_{ic}^t$ has its own interest and merits an independent
emphasis:

\bigskip
Theorem C.\ \ $P_{ic}^t$ is a complex manifold of (complex) dimension
$4g-3$.

\bigskip
   In section 2 we prove some basic results about spherical CR circle
bundles arising from holomorphic line bundles. In section 3 we prove
Theorem C.  We give a representation of the tangent space of $P_{ic}^t$
in the "classical gauge" (see (3.32)), which maps onto the space of
holomorphic (1,0)-forms through $\bar\pa$-operator with the kernel equal to the
space of holomorphic quadratic differentials relative to
the reference Riemann surface. To parametrize our moduli space of spherical
CR structures we introduce a certain local "supporting" manifold in section 4.
We also show the properness of the contact action in our case. In section 5 we
parametrize $\frak  S^t_0$ as a smooth manifold with the aid of the map
$\tau^t$ and the local "supporting" manifold. Finally we prove Theorem A
and Corollary B. On the way to showing Theorem A, we actually obtain another
representation of the tangent space of $P_{ic}^t$, which
is a fourth-order differential equation. (It is basically because the deformation
tensor of spherical CR structures in dimension 3 is of fourth order.) In
Appendix A we prove the $U(1)$-invariant version of Gray's theorem (Theorem
5.1). In Appendix B we give a description of an infinitesimal slice of
$\widetilde{\frak M}_{-1,0}/C_{\hat H}.$

   Our theory for the universal Picard variety has its counterpart in the
Teichmuller theory as shown in the following table:
$$
\centerline{ \vbox { \halign{\strut # & \vrule # \tabskip =5 true pt & # \cr
 \text{Teichmuller space} &  & \text{universal Picard variety}\cr
\noalign{\hrule}
\text{conformal classes} & &\text{spherical CR circle bundles}\cr
\noalign{\hrule}
\text{Riemannian hyperbolic metrics} & & \text{pseudohermitian
hyperbolic geometries}\cr}}}
$$

   Our description of $P_{ic_0}^t$ using $\frak T_{-1,0}$ (combining Theorem
A and Corollary B) has a topological implication. Namely, the topology of
(contact, hence) diffeomorphism group of $\hat S$ in principle can be determined by
the topologies of ${\frak T}_{-1,0}$ and the unimodulo representative space
$\widetilde {\frak M}^0_{-1,0}$. But the topology of ${\frak T}_{-1,0}$ is the
same as that of $P_{ic_0}^t$ by  our theorems, which is well known. To study
the topology of $\widetilde {\frak M}^0_{-1,0}$,
we might define a certain kind of Dirichlet's energy on it and use this
energy functional as our Morse function. A similar strategy works
successfully in studying the topology of Teichmuller space.([Tr])

   Another problem is the analogue of the so-called Nielsen realization
problem about the mapping class group of a Riemann surface. The Nielsen
realization problem says whether any finite subgroup of the mapping class
group $Diff_+/Diff_0$ (of a surface with genus $> 1$) can be "realized"
as a subgroup of $Diff_+$. There is an analytic proof using the above-
mentioned Dirichlet's energy and the so-called Weil-Petersson metric.([Tr])
We wonder if we can do the similar thing for a 3-manifold of circle bundle
type through the study of ${\frak T}_{-1,0}.$

   As we know, the moduli space of Riemann surfaces of fixed genus is the 
quotient of Teichmuller space by the mapping class group.  The 
compactification of the moduli space has been well studied. 
There are a couple of ways to do it. The way using algebraic geometry
was first done by Deligne and Mumford. It was realized later a different approach  
which is based on the Riemannian hyperbolic geometry. (e.g., [SS],
[Pa]) As for the compactification of the universal Picard variety, 
algebraic approaches have been taken up by several authors.
([Ds], [Is], [OS], [Cap], etc.) Towards the problem of 
compactification we hope that
along with the framework of this paper 
there will be a differential-geometric approach 
in the near future.

\bigskip
\n
2. Spherical CR circle bundles

   Let $L$ be a negative holomorphic line bundle over a closed Riemann surface
$N$ of genus $g>1$. For such $N$, there always exists a unique hyperbolic
metric $ds^2$ (i.e. the associated Gaussian curvature $K_{ds^2}=-1$) in its
associated conformal class. Denote the volume form of $ds^2$ by
$\omega_{ds^2}$.  By the Gauss-Bonnet theorem the integral of $\hat
\omega_{ds^2} = -\omega_{ds^2}/2\pi \chi(N) = \omega_{ds^2}/4\pi(g-1)$ equals
1. Hence $[\hat \omega_{ds^2}]$ is a generator  of $H^2(N,Z)=Z$. Write
$c_1(L)=-m[\hat \omega_{ds^2}]$ for $m$ being a positive integer.

\bigskip
\n
Proposition 2.1. There exists a unique (up to a positive constant multiple)
hermitian metric $\|\ \ \|: L\to \Bbb R^+ \cup\{0\}$ such that if we write
$h(z,\bar z)=\|s(z)\|^2$ for a local holomorphic section s, then
$$
i\pa_z\pa_{\bar z}\log h(z,\bar z)=(m/(2g-2))\omega_{ds^2}          \tag 2.1
$$

\bigskip
\n
Proof. Take an arbitrary hermitian metric $\|\ \ \|_0$ and write $h_0(z,\bar
z)=\|s(z)\|^2_0$. Any other hermitian metric $\|\ \ \|^2$ is equal to
$\lambda$$\|\ \ \|^2_0$ for $\lambda$ being a global positive function
defined on $N$. It suffices to solve $\lambda$ for the following equation:
$$
 i\pa_z\pa_{\bar z}\log\lambda(z,\bar z)=(m/(2g-2))\omega_{ds^2} -
                   i\pa_z\pa_{\bar z}\log h_0(z,\bar z)      \tag 2.2
$$

\n
Equating and then multiplying coefficients of $idz\land d\bar z$ in (2.2) by
$g^{1\bar 1}= (g_{1\bar 1})^{-1}$ where $ds^2 = g_{1\bar 1}dzd\bar z$
gives
$$
 \f 12{\Delta}_{ds^2}\log\lambda=\Sigma   \tag 2.3
$$
\n
where $\Sigma$ is a global real function. Note that both
$(m/(2g-2))\omega_{ds^2}$ and $i\pa_z\pa_{\bar z}\log h_0(z,\bar z)$ represent $-2\pi c_1(L)$. It follows that
$\int\Sigma \omega_{ds^2}=0$. So we can solve (2.3) for $\lambda$ unique up
to a positive constant multiple (see, for instance, p.104 in [Au]) and hence
(2.2).  \qquad\qed

\bigskip
   Define the circle bundle $S_L \subset L$ by $\|\ \ \|=1$. The contact bundle
and the CR structure on $S_L$, induced from L, are denoted $H_L, J_L $
respectively. Define the contact form $\theta_L$ on $S_L$ by
$$
\theta_L= -i{\kappa}\pa_L(\|\ \ \|^2)|_{S_L}    \tag 2.4
$$

\n
with the normalizing constant $\kappa=2(g-1)/m$. Locally write $\|ws(z)\|^2=
h(z,\bar z)|w|^2$ for $w$ in $\Bbb C$, a fibre coordinate. A direct computation
using (2.1) shows that
$$
d\theta_L = \pi^*\omega_{ds^2}\ \   (\pi: S_L \subset L\to N\ \ \text{being the
bundle projection})
$$

Let $w^1, w^2$ be orthonormal coframe fields relative to $ds^2$. Let $\theta^1
=w^1 + iw^2,\ \theta^{\bar 1} = w^1-iw^2$ be the corresponding unitary coframe
fields.  Hence $\omega_{ds^2} = w^1\land w^2 = (1/2)i\theta^1\land \theta^{\bar
1}$. From the formulas  on pp. 266-267 in [We2], the pseudohermitian scalar (or
Tanaka-Webster) curvature
$$
\leqalignno{
R_{J_L,\theta_L} & = H^{1\bar 1}R_1{}^1{}_{1\bar 1} \ (\text{since}\
h_{1\bar 1}=\f 12, R_1{}^1{}_{1\bar 1}=(1/2)K_{ds^2})   &(2.5)\cr
                 &= 2(1/2)K_{ds^2} = K_{ds^2} = -1    \cr}
$$

\n
and the torsion
$$
 A_{J_L,\theta_L}=0    \tag 2.6
$$

\n
Therefore by (2.4) in [CL1] the Cartan (curvature) tensor
$$
 Q_{J_L} = 0     \tag 2.7
$$

\n
It follows ([Ca],[CM]) that $(S_L, H_L, J_L)$ is spherical, i.e., locally
CR-equivalent to the unit sphere $S^3$ in $\Bbb C^2$. The following Proposition
shows the uniqueness of the contact form $\theta_L$ in (2.5).

\bigskip
\n
Proposition 2.2. Let $(M, H)$ be a closed contact manifold of dimension 2n+1.
Let $(J, \theta_j), j=1,2$ be pseudohermitian structures on $(M, H)$. (i.e. J
is compatible with H and $\theta_j$'s are contact forms relative to $H$)
Suppose the pseudohermitian scalar curvature  $R_{J,\theta_j} = -1,\ j=1,2$.
Then $\theta_1 = \theta_2$.

\bigskip
\n
Proof. Write $\theta_2 = u^{2/n}\cdot \theta_1$ for $u > 0.\ \ R_{J,\theta_1}$
and $R_{J,\theta_2}$ are related in the following equation:
$$
(2(n+1)/n)\Delta_bu + R_{J,\theta_1}u - R_{J,\theta_2}u^{(n+2)/n}=0  \tag 2.8
$$

\n
(see [JL]; note that $\Delta_b$ is the "negative" sublaplacian relative to
($J, \theta_1$)) Substituting $R_{J,\theta_j} = -1$ in (2.8) gives
$$
\f{2(n+1)}n\Delta_bu = u - u^{(n+2)/n}   \tag 2.9
$$

\n
Suppose u achieves its maximal value $> 1$ at a point $p$. Then we evaluate
(2.9) at $p$:
$$
       0\le \f{2(n+1)}n\Delta_bu = u - u^{\f{n+2}n} < 0
$$

\n
to reach a contradiction. Similarly u cannot achieve its minimal value $< 1$.
Therefore u must be identically equal to 1. \qquad\qed

\bigskip
\n
Corollary 2.3. The map $\iota$ in Corollary B (assuming it is well defined)
is injective.

\bigskip
   Next suppose we have a holomorphic bundle isomorphism $(\phi, f): (L_1, N_1)
\to (L_2, N_2)$ for $[(L_j, N_j)] \in P_{ic},\ j=1,2$. By Proposition
2.1 and noting that the biholomorphism $f: N_1\to N_2$ is an isometry, we
conclude that
$$
\phi^*(\|\ \ \|_{L_2}) = c\|\ \ \|_{L_1}    \tag 2.10
$$

\n
for some constant $c > 0$, where $\|\ \ \|_{L_j},\ j=1,2$ denote hermitian
metrics obtained in Proposition 2.1 with respect to $L_j$. Let $m_c$ denote the
multiplication by $c$ on the line bundle. Thus by (2.10) the composition
$\phi \circ m_c^{-1}: S_{L_1}\to S_{L_2}$ is a CR equivalence. We have shown
that the map $\tau:[(L,N)]\to [(S_L,H_L,J_L)]$ is well-defined. Furthermore, we
have

\bigskip
\n
Proposition 2.4. The map $\tau$ is injective.

\bigskip
\n
Proof. Let $D_L$ denote the disc bundle of $L$ with the boundary $\pa D_L
= S_L$.  Since $S_L$ is the strictly pseudoconvex boundary of the complex
manifold  $D_L$, it is CR-embeddable in $\Bbb C^N$ and coordinate functions (CR
functions on $S_L$) extend holomorphically to $D_L$. (see, e.g., Theorem 5.3 in
[K1], Corollary of Theorem 1.3 in [K2], p.91(5.3.5) in [FK]) So we have a map
$\psi: \bar D_L \to \Bbb C^N$, holomorphic in $D_L$, and CR equivalent between
$\pa \bar D_L =  S_L$  and $\psi(S_L)$. Denote the 0-section of $D_L$ by
$\Sigma.\ \Sigma$ is biholomorphic to the closed Riemann surface $N$. We claim
$\psi:D_L\backslash\Sigma\to\Bbb C^N$ is a biholomorphism onto its image. First
observe that $\psi$ is biholomorphic near the boundary $S_L$ and the disc
bundle $D_L(\rho) = \{s \in D_L: \|s\|_L <  \rho\}$ of radius $\rho, 0< \rho <
1$, is strictly pseudoconvex. By continuity,  there exists a smallest $\rho_0
\ge 0$ such that $\psi$ is biholomorphic on  $D_L\backslash D_L(\rho_0)$ and
fails to be biholomorphic on $S_L(\rho_0) = \pa\overline{D_L(\rho_0)}$. Suppose
$\rho_0 > 0$. Take $q \in S_L(\rho_0)$. Near q consider the determinant of the
Jacobian matrix of $\psi$, denoted $J_\psi$. If $J_\psi(q) = 0$, the subvariety
defined by  $J_\psi = 0$ must contain a point near $q$ but out of
$\overline{D_L(\rho_0)}$ due to pseudoconvexity of $D_L(\rho_0)$, which
contradicts $\psi$ being biholomorphic on
$D_L\backslash\overline{D_L(\rho_0)}$
(where $J_\psi \ne 0)$. Thus $J_\psi(q)\ne 0$. Hence $\psi$ is a local
biholomorphism near $S_L(\rho_0)$. Therefore $\psi$ must be globally injective
near $S_L(\rho_0)$ since it is biholomorphic on "one side" of $S_L(\rho_0)$. In
conclusion $\rho_0$ must be 0 and we have proved our claim.

   Now take two holomorphic line bundles $(L_j, N_j),\ j=1,2$ with associated
spherical CR circle bundles $S_{L_j}$'s being isomorphic. That is to say, there
exists a CR equivalence $\phi: S_{L_1}\to S_{L_2}$. As just discussed above,
there exists a map $\psi_1: \bar D_{L_1}\to \Bbb C^N$, biholomorphic on
$D_{L_1}\backslash\Sigma_1\ (\Sigma_j$'s, $j=1,2$, denote zero sections of
$L_j$ respectively) Moreover the CR embedding $\psi_1 \circ \phi^{-1}:
S_{L_2}\to  \Bbb C^N$ extends to a map $\psi_2: \bar D_{L_2} \to \Bbb C^N$,
biholomorphic on $D_{L_2} \backslash\Sigma_2$, with Range $\psi_2 =$ Range
$\psi_1$ by the uniqueness of solution for the complex Plateau problem in
$\Bbb C^N$. ([HL], [Y]) Since $\Sigma_j$'s are biholomorphic to closed Riemann
surfaces $N_j$'s respectively, $\psi_j(\Sigma_j)$ consists of a point $p_j$ in
$\Bbb C^N$. Suppose $p_1 \ne p_2$. Take  a suitable neighborhood $U$ of $p_2$
such that $\psi_1^{-1}(U\backslash p_2)$ is biholomorphic to $\psi_2^{-1}
(U\backslash p_2)$.  But they  have different topological types since the
latter  is a tubular neighborhood of  a closed Riemann surface $\Sigma_2$. So
$p_1 = p_2$ and $\phi_1 = \psi_2^{-1} \circ \psi_1 : D_{L_1}
\backslash\Sigma_1\to D_{L_2}\backslash\Sigma_2$ is a
biholomorphism.  Furthermore it is easy to see that a punctured fibre disc must
be mapped by $\phi_1$ onto a punctured disc with the puncture sitting in
$\Sigma_2$. (just noting that the puncture in $\Sigma_1$ is a removable
singularity) We therefore extend $\phi_1$ to a map (still denoted $\phi_1$)
from $D_{L_1}$ into $D_{L_2}$ carrying $\Sigma_1$ into $\Sigma_2$. We claim
$\phi_1$ is continuous on $\Sigma_1$. Take $q \in \Sigma_1$ and $\widetilde q =
\phi_1(q) \in \Sigma_2$. Centered at $q$, we have local holomorphic coordinates
$(z, w) \in D \times D^*$ for $\Sigma_1 \times$ fibres where $D\ (D^*$ resp.)
denotes the (punctured resp.) disc. Given neighborhoods $U$ and $V$ of
$\widetilde q$ in $D_{L_2}$ with $\bar V \subset U$, there exists a positive
number $r$ such that $\{(0,w): |w|<r\}$ is mapped into $V$. Observe that
$$
d_{D_{L_2}\backslash\Sigma_2}(\phi_1(z,w),\phi_1(0,w))\le
                   d_{D\times D^*}((z,w),(0,w))\le d_D(z,0)    \tag 2.11
$$

\n
where "$d$" denotes the Kobayashi distance. Let $b =
d_{D_{L_2}\backslash\Sigma_2}((D_{L_2}\backslash\Sigma_2)\cap
V,(D_{L_2}\backslash \Sigma_2)\backslash U) > 0$.
Then there is a positive number $r'$ such that $d_D(z, 0) < b$ for $|z| <
r'$,  so it follows by (2.11) that $\phi_1(z, w)$ is in $U$ for $|z| < r'$,
$|w| < r$. Once we know $\phi_1$ is continuous on $\Sigma_1$, then it must be
holomorphic on $D_{L_1}$ by the Riemann extension theorem. (we can also just
invoke Theorem 6.2 in [Ko] p.93 to replace the above argument) Similarly extend
$\phi_1^{-1}$ holomorphically to $D_{L_2}$. Since the holomorphic map
$\phi_1^{-1}\circ \phi_1=$ identity on $D_{L_1}\backslash\Sigma_1$, it must be
an identity on $D_{L_1}$. We have shown that $\phi_1$ is
a biholomorphism between $D_{L_1}$ and $D_{L_2}$. Define biholomorphisms
$\phi_\rho : D_{L_1}(\rho)\to D_{L_2}(\rho)$ for $\rho > 0$ by $\phi_\rho(y) =
\rho \phi_1(y/\rho)$.  In local coordinates $(z, w)$ with w being fibre
coordinate, we can write $\phi_\rho = (\widetilde z, \widetilde w)$ as a
function of $y=(z, w)$: at $z=0$,
$$
\widetilde z = O(w/\rho),   \widetilde w = cw + O(w^2/\rho)
$$

\n
for some nonzero constant $c$. As $\rho$ goes to infinity, $\widetilde z$
approaches to 0 and $\widetilde w$ goes to $cw$. That is to say,
$\lim_{\rho\to\infty}\phi_\rho = \phi_\infty : L_1\to L_2$
exists and is a linear isomorphism on each fibre, and from the above
argument $\{\phi_\rho=(\widetilde z, \widetilde w)\}$ is uniformly bounded on
any compact coordinate neighborhood around a point. It follows that
$\phi_\infty$ is holomorphic. Apply
a similar argument to $\phi_\rho^{-1}(x)=\rho \phi_1^{-1}(x/\rho)$. We obtain a
holomorphic map $\psi = \lim \phi_\rho^{-1} : L_2\to L_1$ and it is easy to see
that $\psi \circ \phi_\infty=\phi_\infty \circ \psi=$ identity. Therefore
$\phi_\infty$ is a holomorphic bundle isomorphism.
\qquad  \qed

\bigskip
\n
3. Parametrizing $P_{ic}^t$ as complex manifold: Proof of Theorem C

   First we describe $P_{ic}$ by complex structures with special
properties on the fixed background line bundle $\hat L$ (considered as a
smooth line bundle). Since every holomorphic line bundle $(L,N)$ of the fixed
$c_1$ is isomorphic to $(\hat L,\hat N)$ as smooth line bundle, complex
structures on $L$ and $N$ are pulled back to $\hat L$ and $\hat N$
respectively. Let $Bdiff$  denote the group of smooth bundle automorphisms of
$(\hat L,\hat N)$. Let $J, c$ denote complex structures on $\hat L, \hat N$
respectively. The space of all $((\hat L,J),(\hat N,c))$
such that the projection from $\hat L$ onto $\hat N$ is holomorphic with
respect to $(J,c)$ modulo $Bdiff$ is in one-to-one correspondence with
$P_{ic}.$   Let $m_\rho:\hat L\to \hat L$  denote the fibre multiplication by
$\rho$, a complex number.  Let $\Bbb C^*$ denote the subgroup of
$Bdiff$, consisting of all $m_\rho$ with nonzero $\rho$. Let  $\hat J$ denote
the complex structure on $\hat L$ (and also on $\hat S$, cf.
section 1) associated to the fixed (or reference) holomorphic line bundle
$(\hat L,\hat N)$. On $\hat L$ we consider the space $\widetilde P_{ic}$ of all
smooth almost  complex structures $J$ respecting the same orientation as given
by $\hat J$ on $\hat L$ and satisfying the following conditions:

\bigskip
\n
(3.1)\ \  $m_\rho^*J =J$ for nonzero $\rho$ (i.e. $J$ is $\Bbb C^*$-invariant)
and

\n
(3.2) on fibres, $J$ is induced by the usual complex structure on $\Bbb C$
      in local trivializations.

\bigskip
\n
Proposition 3.1. Any $J$ in $\widetilde P_{ic}$ is integrable.

\bigskip
\n
Proof. First fix a system of local coordinates $(z,w)$ on $\hat L$ (holomorphic
with respect to the original reference complex structure $\hat J$ on $\hat L$)
with fibre coordinate $w$. Let $\pa_z, \pa_w$ denote the tangent vectors
$\pa/\pa z, \pa/\pa w$ respectively for short. The condition (3.1) allows us to
construct a $\Bbb C^*$-invariant (1,0)(relative to $J$) tangent vector $Z_1$ on
$\hat L$ by moving a chosen (1,0) section by the action of $\Bbb C^*$. Write
$$
Z_1 = f\pa_z + g\pa_{\bar z} + hw\pa_w + \ell\bar w\pa_{\bar w}.
$$

\n
$\Bbb C^*$-invariance implies that $f, g, h, l$ are smooth functions only in
$z, \bar z$.  It follows that $[Z_1, \pa_w] = -h\pa_w$. Now we can compute the
Nijenhuis tensor:
$$
N(Z_1, \pa_w) = -4[Z_1, \pa_w] -4iJ[Z_1, \pa_w] = 0,
$$

\n
and it is easy to see that
$$
N(\pa_w, \pa_w) = N(\pa_w, \pa_{\bar w}) = N(Z_1, \pa_{\bar w}) = N(Z_1,Z_{\bar
1}) = 0.
$$

\n
Thus by noting that $N$ is skew-symmetric, the Nijenhuis tensor vanishes.
\qquad\qed

\bigskip
   Observe that $(\hat L,J)$ in $\widetilde P_{ic}$ can be pushed down to
$(\hat N,c)$ for some complex structure $c$. (for $v$ tangent to $\hat N, c(v)$
is defined to be $\hat\pi_*J(s_{0*}(v))$.  Here $s_0$ is the zero section of
$\hat L$ over $\hat N$ and $\hat \pi: \hat L\to \hat N$ is the natural
projection. It follows that $\hat \pi$ is holomorphic with
respect to $(J,c))$.  Hence the quotient space $\widetilde P_{ic}/Bdiff$ is in
one-to-one correspondence with $P_{ic}$. Let $Bdiff_0$ denote the group
of smooth bundle automorphisms $(\widetilde\phi,\phi)$ of $(\hat L,\hat N)$
with $\phi: \hat N\to \hat N$ being isotopic to
the identity. Denote the quotient group $Bdiff_0/\Bbb C^*$ by $\frak B$. (note
that $\Bbb C^*$ is contained in the center of $Bdiff_0$) Define the
Teichmuller-type space $P_{ic}^t$ to be $\widetilde P_{ic}/{\frak B} =
\widetilde P_{ic}/Bdiff_0$. We are going to show that $\frak B$ acts freely and
properly on $\widetilde P_{ic}$ and $P_{ic}^t$ can be parametrized as complex
manifold. First parametrize $\widetilde P_{ic}$ and $\frak B$. Let us do a
priori computation of the tangent space of $\widetilde P_{ic}$
at a reference point $\hat J$. Denote $J_t$ a family of elements in
$\widetilde P_{ic}$ with $J_0 = \hat J$. Let $E$ be the derivative of $J_t$ in
$t$ at $t=0$ (considered in the space of
endomorphisms of $T\hat L$). $J_t$ being almost complex structures implies that
$E$ satisfies the following equation:
$$
E\circ \hat J + \hat J \circ E = 0    \tag 3.3
$$

\n
Take local holomorphic coordinates $(z, w)$ relative to $\hat J$ as in the proof of
Proposition 3.1. Write $\pa_b = \pa/\pa z^b, b=1,2$ for short where $z^1 = z,
z^2 = w$. We express $E$ as below:
$$
E = \Sigma E_a{}^b dz^a \otimes \pa_b + E_a{}^{\bar b}dz^a \otimes
       \pa_{\bar b} + \text{conjugate}.
 $$

\n
It follows from (3.3) that
$$
E_a{}^b = 0 (hence E_{\bar a}{}^{\bar b} = \overline{(E_a{}^b)} = 0)   \tag 3.4
$$

\n
Condition (3.2) means $J_t(\pa_2) = i\pa_2$ whose differentiation in $t$ at
$t=0$ gives $E(\pa_2$ or $\pa_{\bar 2}) = 0$. Hence we have
$$
E_2{}^{\bar 1} = E_2{}^{\bar 2} = 0   \tag 3.5
$$

\n
Besides, differentiating (3.1) tells that $E$ is $\Bbb C^*$-invariant.
Therefore both $E_1{}^{\bar 1}$ and $E_1{}^{\bar 2}/\bar w$ are independent of
the variable $w = z^2$. Together with (3.4),(3.5) we obtain
$$
E = E_1{}^{\bar 1}(z,\bar z)dz\otimes \pa_{\bar z} + E_1{}^*(z,\bar z)\bar w
dz\otimes \pa_{\bar w} + \text{conjugate}     \tag 3.6
$$

\n
where $E_1{}^*(z,\bar z)$ is just $E_1{}^{\bar 2}/\bar w. \ \ E_1{}^{\bar 1}$
and $E_1{}^*$ satisfy the transformation law:
$$
\widetilde E_1{}^{\bar 1}= E_1{}^{\bar 1}\overline{(h')}(h')^{-1}   \tag 3.7.1
$$
$$
\widetilde E_1{}^* = E_1{}^*(h')^{-1} + E_1{}^{\bar
1}(h')^{-1}\overline{g'g^{-1}}   \tag 3.7.2
 $$

\n
under the coordinate change of trivializations:
$$
\widetilde z = h(z), \widetilde w = g(z)w    \tag 3.8
$$

\n
for biholomorphic $h$ and nonzero holomorphic $g$. Therefore we can talk about
smooth or $H^s$ (Sobolev $s$-norm bounded) $E$ if $E_1{}^{\bar 1}$ and
$E_1{}^*$ are smooth or $H^s$. (the Sobolev $s$-norm can be defined via either
a chosen partition of unity or a chosen covariant derivative on $\hat N$)
Similarly by conditions (3.1), (3.2), we can write an element $J$ in
$\widetilde P_{ic}$ as
$$
J = \Sigma J_1{}^bdz\otimes \pa_b + idw\otimes \pa_w +
             \text{conjugate}
$$

\n
where b runs over $1,\bar 1,2,\bar 2$ and $J_1{}^1, J_1{}^{\bar 1}, J_1{}^2/w,
J_1{}^{\bar 2}/\bar w$ are independent of $w$. Therefore we can talk about $H^s
J$ if these components are all in $H^s$. Let $\widetilde P_{ic_s}$ denote the
set of all such $H^s J$. Let $\frak E_{\hat J} (\frak E^s_{\hat J}$ resp.)
denote the linear space of all smooth ($H^s$, resp.) tensors $E$ of the type
(3.6). Since $\hat N$ is compact, ${\frak E}_{\hat J}$ is a tame Frechet space
in the terminology of [H] while $\frak E^s_{\hat J}$ is a Hilbert (hence
Banach) space. Define a map $\Phi_{\hat J}: \frak E_{\hat J} \to\widetilde
P_{ic}$ by $$
\Phi_{\hat J}(E) = (I - (1/2) E \circ \hat J)\circ\hat J\circ(I - (1/2)E\circ
\hat J)^{-1}     \tag 3.9
$$

\n
for small (in $C^\infty$-topology) $E$. It is easy to see that $\Phi_{\hat J}$
extends to $\frak E^s_{\hat J}$ (still denoted $\Phi_{\hat J}$) with the range
$\widetilde P_{ic_s}$ for large enough s, say, $s \ge 2$ by the Sobolev lemma.
(for $s\ge 2, H^s$-space is contained in $C^0$ and forms an algebra. Note also
that the inverse of a nonzero $H^s$ function on $\hat N$ is still in $H^s$)
Moreover $\Phi_{\hat J}$ is injective for small $E$ in $\frak E^s_{\hat
J}$ as the inverse $\Phi^{-1}_{\hat J}$ can be given precisely by
$$
\Phi^{-1}_{\hat J}(J) = 2(J -\hat J)(J +\hat J)^{-1}\hat J,    \tag 3.10
$$

\n
and it is easy to compute that $(d/dt)\Phi_{\hat J}(tE) = E$ at $t=0$.
(consider $\widetilde P_{ic}$ sitting in $End(T\hat L$)) We use $\{\Phi_J: J
\in \widetilde P_{ic}$ or $\widetilde P_{ic_s}\}$ to parametrize $\widetilde
P_{ic}$ or $\widetilde P_{ic_s}$. The transition map for the overlap region
have the precise formula by composing (3.10) and (3.9) for two different $\hat
J$'s. Observe that, with respect to a basis, each component of the transition
map is a polynomial in components of $E$. It follows that the transition map
is $C^\infty$ (smooth tame in the smooth case) and hence a $
C^\infty$-diffeomorphism by symmetry. We have proved

\bigskip
\n
Proposition 3.2. $\widetilde P_{ic_s} (\widetilde P_{ic}$, resp.) is a smooth
Hilbert (tame Frechet, resp.) manifold for $s \ge 2$.

\bigskip
   Next a priori computation shows that a tangent vector of $Bdiff_0$
at the identity has the following form:
$$
X = v^1 \pa_z + v^{\bar 1} \pa_{\bar z} + v^* w \pa_w + \bar{v^*}\bar w
\pa_{\bar w}   \tag 3.11
$$

\n
in a local trivialization $(z,w)$ as above, where $v^1$ and $v^*$ are
independent of $w$ and satisfy the following transformation law:
$$
\cases \widetilde v^1 = v^1 h'(z)\\
       \widetilde v^* = v^* + v^1 g'(z)g(z)^{-1}\endcases \tag 3.12
$$

\n
for the change of trivializations (3.8). Let $\widetilde{\frak V}^s$ denote
the Hilbert space of all $X$ satisfying (3.12) with bounded $H^s$-norm.
(may be defined by fixing a finite number of trivializations and a
corresponding partition of unity for $\hat N$ so that the $H^s$-norm is locally
provided by the sum of $H^s$-norms of $v^1$ and $v^*$)  On the other hand a
bundle automorphism $\phi$ of $\hat L$ can be expressed as
$$
\phi: (z,w)\to(\psi(z,\bar z), \lambda(z,\bar z)w)
$$

\n
in trivialization $(z,w)$, where $\psi, \lambda$ obey the following
transformation law:
$$
\align
&  \widetilde\psi(h(z),\overline{h(z)}) = h(\psi(z,\bar z))\\
&  \widetilde\lambda(h(z),\overline{h(z)}) = {\lambda}(z,\bar z)g(\psi(z,\bar
z))g(z)^{-1} \endalign
$$

\n
according to (3.8). We say $\phi$ is $H^s$ if $\psi$ and $\lambda$ are $H^s$
for each trivialization. Let $Bdiff_0^s$ denote the topological space of all
$H^s$ bundle automorphisms of $\hat L$ with obvious $H^s$-topology.
Take $X$ in $\widetilde{\frak V}^s$. We want to associate a bundle automorphism
$\phi_X$ in $Bdiff_0^s$. Take a (smooth) metric $ds^2$ on $\hat N$ and a
hermitian connection $\nabla$ of $\hat L$ over $\hat N$. Let $V=\hat\pi_*X$
be the projection of $X$ on $\hat N$. Locally $V = v^1 \pa_z + v^{\bar 1}
\pa_{\bar z}$ if $X$ is
expressed as in (3.11). Let $\gamma(p, V(p), t)$ be the geodesic relative to
$ds^2$ with initial point $p$ and initial velocity $V(p)$. It is well known
that
$\gamma$ is smooth in $(p,v,t)$ where defined. Let $s_0$ denote the local section
of $\hat L$ given by $(z, 1)$. Define the connection form $\Gamma$ on $\hat N$ by
$$
{\nabla}_v(s_0) = \Gamma(v)s_0    \tag 3.13
$$

\n
for tangent vectors $v$ on $\hat N$. Denote $(d/dt){\gamma}(p,V(p),t)$ by
$\gamma'_t(p,V(p),t)$. In trivialization $(z,w)$, we identify $p$ with $z_0$
and move the fibre element $w_0$ parallelly according to (3.13) along the
geodesic $\gamma(z_0,V(z_0),t)$ (instead of
$\gamma((z_0,\bar z_0),V(z_0,\bar z_0),t)$
for short) to get $w_1$ at time = 1 (for small $V$). It is then easy to compute
$$
 w_1 = w_0 \exp[-\int^1_0 \Gamma(\gamma'_t(z_0,V(z_0),t))dt]    \tag 3.14
$$

\n
(3.14) suggests the following choice of $\phi_X$:
$$
 \phi_X(z_0,w_0)=(\gamma(z_0,V(z_0),1),   \
                       w_1 \exp[v^*(z_0) + \Gamma(V(z_0))]).   \tag 3.15
$$

\n
Here we write $v^*(z_0)$ instead of $v^*(z_0,\bar z_0)$ and recall
$V=v^1\pa_z+v^{\bar 1} \pa_{\bar z}$ and $v^*$ are local components of
$X$ as expressed in (3.11). We claim the definition of $\phi_X$ given by (3.15)
is independent of the choice of trivialization. Let $(\widetilde z,\widetilde
w)$ be another trivialization related to $(z,w)$ by (3.8). We have
corresponding local section $\widetilde s_0$ given by $(\widetilde z,1)$ and
associated connection form $\widetilde\Gamma$. It is easy to see that
$g(z)\widetilde s_0 = s_0$ and
$$
\widetilde\Gamma = \Gamma - dg\cdot g^{-1}.     \tag 3.16
$$

\n
Now applying (3.16) to $\gamma'_t(z_0,V(z_0),t)$ gives $\widetilde w_1 =
g(z_1)w_1$ where $z_1 = \gamma(z_0,V(z_0),1)$. (note that $\widetilde w_0
=g(z_0)w_0$ and $\widetilde w_1$ is given according to (3.14)).  From the
transformation law (3.12)  for $v^*$ we can easily show that $v^*(z_0) +
\Gamma(V(z_0))$ is invariant under the change of trivialization (3.8).
Altogether we have proved our claim. Observe that $\Gamma$ is smooth and
$\gamma$, hence $\gamma'_t$ is also smooth in their arguments. It follows that
${\phi}_X$ is $H^s$ if $X$ is $H^s$. So we have defined a map $\Sigma:
{\widetilde{\frak V}}^s\to Bdiff_0^s$ by $\Sigma(X)={\phi}_X$. If we write
$$
\phi_X(z,w)=({\phi}^1_X(z,\bar z), w \exp{\phi}^2_X(z,\bar z)),
$$

\n
then ${\phi}^1_X$ gives rise to a global diffeomorphism on $\hat N$ (still
denoted ${\phi}^1_X)$ and the inverse of $\Sigma$ can be given by
$$
\leqalignno{& V=  P^{-1}({\phi}^1_X) &(3.17)\cr
            & v^* = \phi^2_X - \Gamma(V(\cdot) + \int^1_0
             \Gamma(\gamma'_t(\cdot,V(\cdot),t))dt. \cr
            &    \text{(with}\ V \ \text{replaced by the first formula)} \cr}
$$

\n
Here $P$ is the usual map of parametrization from vector fields to
diffeomorphisms on $\hat N$ via the geodesic flow. Now it is clear that
$\Sigma$ is a homeomorphism from an open set of small $X$ to a neighborhood
of the identity, say, $U$. Let $l_\psi$ denote the composition with $\psi$ from
the left. By composing $\Sigma$ with $l_\psi$ for smooth elements $\psi$ in
$\frak B^s$, we obtain an atlas $\{(l_\psi(U), \Sigma^{-1} \circ
l_{\psi^{-1}}): \psi$ is a smooth element in $Bdiff_0^s\}$ for $Bdiff_0^s$.
(note that the set of smooth elements is dense in $Bdiff_0^s$, and the
composition map and the map taking each diffeomorphism to its inverse are
$C^0$) To show the transition map being smooth is a matter of direct
verification (using (3.14), (3.15), (3.17)) : one only has to observe that
composing with a smooth element is smooth in the original argument. (Actually
we can prove $Bdiff_0^s$ is a topological group for $s\ge 3$ under the
operation of composition of $H^s$-maps (cf. section 3 in [Eb]). However the
composition map is only $C^0$ but not $C^1$, so to get $C^\infty$
differentiable structure on $Bdiff_0^s$, we have to restrict to smooth
elements as our "centers" of charts) Let $\frak V^s$
denote the quotient space $\widetilde{\frak V}^s/C$ where $C$ consists of all
$X$ in (3.11) with $v^1=0, v^*=a$ constant complex number.(this is well defined
according  to the transformation law (3.12)) Since any finite-dimensional
subspace  of a Hilbert space is closed, we can identify $\frak V^s$ with the
closed orthogonal complement of $C$ in $\widetilde{\frak V}^s$, which
inherits the  Hilbert space structure from $\widetilde{\frak V}^s$. Recall
$\frak B^s = Bdiff_0^s/\Bbb C^*$ where $\Bbb C^*$ consists of all fibre
dilations by
nonzero constant complex numbers.  (see the beginning of this section) Observe
that $C$ is mapped into $\Bbb C^*$ by $\Sigma$ through the exponential function
according to (3.14), (3.15), so $\Sigma$ induces a homeomorphism from a
neighborhood of 0 in $\frak V^s$ to a neighborhood of the identity in
$\frak B^s$. Similar construction as for $Bdiff_0^s$ gives us the desired
charts for $\frak B^s$. We have proved

\bigskip
\n
Proposition 3.3. $Bdiff_0^s$ and $\frak B^s$ are smooth Hilbert manifolds.

\bigskip
   Next we consider the behavior of $\frak B^{s+1}\ \ (\frak B$, resp.)
acting on $\widetilde P_{ic_s}\ (\widetilde P_{ic}$, resp.) by the pullback.
(well-defined because $C^*$ is contained in the center of $Bdiff_0$) First we
have

\bigskip
\n
Proposition 3.4. ${\frak B}^{s+1}$ acts freely on $\widetilde P_{ic_s}$ for
$s\ge 4$; in particular,  $\frak B$ acts freely on $\widetilde P_{ic}$.

\bigskip
\n
Proof. A bundle automorphism $\phi$ in $H^{s+1}$ fixing a complex  structure
in $\widetilde P_{ic_s}$ can be pushed down to an $H^{s+1}$-biholomorphism
$\underline\phi$ relative to the pushed down $H^s$-complex structure on $\hat
N$. Since $\underline\phi$ is isotopic to the identity map, it follows by a
standard result for genus $\ge 2$ Riemann surfaces (e.g. p.39 in [Tr]) that
$\underline\phi$ must itself
be the identity map.  (to apply the Newlander-Nirenberg theorem [NN, Theorem
1.1], we require $s \ge 4$ so that $H^s$ is contained in $C^2$)  Thus $\phi$
is just a fibre multiplication by a nonzero holomorphic function $\lambda$ on
$\hat N$. Compactness of $\hat N$ implies $\lambda$ must be a constant $\rho$.
Therefore $\phi = m_\rho$ belongs to $C^*.$ \qquad \qed

\bigskip
\n
Proposition 3.5. $\frak B^{s+1}$ acts properly on $\widetilde P_{ic_s}$ for
$s\ge 4$ : i.e. if $\phi_j^*J_j = \widetilde J_j$ converges to $\widetilde
J$ and $J_j$ converges to $J$ in $H^s$  with $J_j$ in $\widetilde P_{ic_s},
[\phi_j]$ in ${\frak B}^{s+1}$, then there exists a  subsequence of $[\phi_j]$
which converges in $H^{s+1}$ to some $[\phi]$.

\bigskip
\n
Proof. First we have ${\phi_j}^*J_j = \widetilde J_j$, so $\phi_j$ can be
pushed down to a biholomorphism $\underline\phi_j$ from $(\hat N, \widetilde
c_j)$ to $(\hat N,c_j)$, where $\widetilde c_j, c_j$ are pushed
down complex structures of $\widetilde J_j,  J_j$ respectively. There is a
diffeomorphism between $H^s$ oriented complex structures and $H^s$ hyperbolic
metrics of Gaussian curvature -1 on a closed surface of genus $\ge 2$ ([Tr]),
so we can apply the Ebin-Palais theorem (Theorem 2.3.1 in [Tr]) to conclude
that a subsequence of $\underline\phi_j$ converges in
$H^{s+1}$ to some $\underline\phi$. Let $\underline\phi(p)=q$ for $p,q$ in
$\hat N$. Take holomorphic coordinates $z, \widetilde z$ with respect
to $\widetilde c= \lim \widetilde c_j, c= \lim c_j$ around p,q respectively so
that $\underline\phi$ satisfies the $\bar\pa$-equation in these coordinates.
Take local trivializations $(z,w)$ and $(z',w')$ of $\hat L$ (which may not be
holomorphic with respect to $\widetilde J$ and $J$ respectively). Write
$\phi_j$ in these local trivializations:
$$
\phi_j: (z, w)\to (z', w')=(\underline\phi_j(z,\bar z),u_j(z,\bar z)w)
$$

\n
for large $j$. Here $\underline\phi_j$ tends to  $\underline\phi$ in
$H^{s+1}$ as $j$ goes to infinity. Moreover since $\underline\phi$ is
holomorphic with respect to $z, \pa_{\bar z}\underline \phi_j$
goes to zero in $H^s$. (note that we need the $H^s$ version of the
Newlander-Nirenberg theorem ([FK]) to
conclude that $\underline \phi_j$ is still in $H^{s+1}$ with respect to
the $z$-coordinate) Now write $J_j$ in $(z',w')$:
$$
J_j = dz' \otimes (f_j \pa_{z'} + g_j \pa_{\bar z'} + h_j w'\pa_{w'} + \ell_j
\bar w'\pa_{\bar w'})+idw' \otimes \pa_{w'} + \text{conjugate}
$$

\n
where $f_j, g_j, h_j, \ell_j$ are $H^s$ functions in $z,\bar z$ according
to (3.1),(3.2).  Moreover $J_j^2 = -I$ implies that $f_j, g_j, h_j,\ell_j $
satisfy the following algebraic conditions:
$$
\leqalignno{\text{(a)}\ \  & f_j^2 + |g_j|^2 = -1   &(3.18)\cr
            \text{(b)}\ \  & g_j(f_j + \bar f_j) = 0  \cr
            \text{(c)}\ \  & h_j(f_j + i) + g_j \bar\ell_j = 0 \cr
            \text{(d)}\ \  & \ell_j(f_j - i) + g_j \bar h_j = 0  \cr}
$$

\n
Similarly for $\widetilde J_j$ we write capital $F_j, G_j, H_j, L_j$ for
corresponding coefficients in trivialization $(z, w)$. Computing $\phi_j^*J_j$
and comparing corresponding coefficients of $wdz \otimes \pa_w$ and $\bar w
dz \otimes \pa_{\bar w}$ with $\widetilde J_j$, we obtain
$$
\leqalignno{\text{(a)}\ \ & (i-f_j+e_j^1)\f\pa{\pa z}
           (u_j)-[\underline\phi'(z)(\overline{\underline \phi'(z))}^{-1}
           g_j +e_j^2]\f\pa{\pa\bar z}(u_j) = u_j \widetilde H_j   &(3.19)\cr
            \text{(b)}\ \ & (-i-f_j+e_j^1)\f\pa{\pa z}(\bar
         u_j)-[\underline\phi'(z)(\overline{\underline\phi'(z)})^{-1}g_j +
           e^2_j]\f\pa{\pa\bar z}(\bar u_j) = \bar u_j\widetilde L_j\cr}
$$

\n
where $\widetilde H_j = H_j - \f{\pa\underline \phi_j}{\pa z}h_j -
\f{\pa\bar{\underline \phi}_j}{\pa
z}\bar\ell_j, \widetilde L_j = L_j - \f{\pa\underline\phi_j}{\pa z}
\ell_j - \f{\pa\bar{\underline \phi_j}}{\pa z}\bar h_j$ and
$$
e_j^1 = \f{\pa z}{\pa \bar z'}\f{\pa\bar z'}{\pa z}(f_j-\bar f_j) -
        \f{\pa z}{\pa \bar z'}\f{\pa z'}{\pa z}g_j-
        \f{\pa z}{\pa z'}\f{\pa\bar z'}{\pa z}\bar g_j
$$
$$
e_j^2 = \f{\pa\bar z}{\pa z'}(\f{\pa z'}{\pa z}f_j + \f{\pa\bar z'}{\pa z}\bar
         g_j) - [\underline\phi'(\bar{\underline\phi'})^{-1} -
        \f{\pa z'}{\pa z}(\overline{\f{\pa z}{\pa z'}})]g_j
    + (\overline{\f{\pa z}{\pa z'}})\f{\pa\bar z'}{\pa z}\bar f_j
$$

\n
Here $\f{\pa z'}{\pa z}\ (\f{\pa z}{\pa z'}$, etc., resp.) means $\f\pa{\pa
z}(\underline\phi_j)(\f\pa{\pa z'}(\underline\phi_j)^{-1}$,
etc., resp.). It is easy to see that $e_j^1$ and $e_j^2$ converge to zero
in $H^s$ as $j$ goes to infinity since $\underline\phi_j$ goes to a
biholomorphism  $\underline\phi$, and obviously $\widetilde H_j$ and
$\widetilde L_j$ converge to $H-(\f{\pa\underline \phi}{\pa z})h$ and
$L-(\f{\pa\underline\phi}{\pa z})\ell$ in $H^s$
respectively, where $H=\lim H_j, h=\lim h_j, L=\lim L_j,\ell=\lim \ell_j$. Let
$$
\align
& D_j = (i - f_j)\f\pa{\pa z} -\underline\phi'(\bar{\underline\phi'})^{-1}
        g_j\f\pa{\pa \bar z}\\
& D'_j= (-i - f_j)\f\pa{\pa z} -\underline\phi'(\bar{\underline\phi'})^{-1}
        g_j\f\pa{\pa \bar z}\endalign
$$

\n
and
$$
\align
& D  = (i - f)\f\pa{\pa z} -\underline\phi'(\bar{\underline\phi'})^{-1}
        g\f\pa{\pa \bar z}\\
& D' = (-i - f)\f\pa{\pa z} -\underline\phi'(\bar{\underline\phi'})^{-1}
        g\f\pa{\pa \bar z}\endalign
$$

\n
where $f = \lim f_j, g = \lim g_j$. Taking the limit of (3.18)(b) gives $g = 0$
or $Re f = 0$. If $g = 0$ at $p, f = i$ or $-i$ by the limiting form of
(3.18)(a). Then either $D$ or $D'$ is not zero at $p$ and equals $\pm
2i\f\pa{\pa z}$. Hence either $D$ or $D'$ is elliptic in a neighborhood of $p$.
On the other hand, if $g$ does not vanish at $p$, then $Re f = 0$ at
$p$. Suppose $(i-f)(a-ib)-\underline\phi'(\bar{\underline\phi'})^{-1}g(a+ib)=0$
at $p$ for real nonzero $a$ or $b$. Then it follows that the absolute value of
$(i-f)/\underline\phi'(\bar{\underline\phi'})^{-1}g$ equals 1. By
(3.18)(a)(limiting version) and $Re f =0$ at $p$, we get $f=i$,
which implies $g=0$ by (3.18)(a) again, a contradiction. Therefore $a=b=0.$
We have proved that $D$ is elliptic around $p$ in the case of $g(p)$ not equal
to zero, so in any case we use either $D$ or $D'$ to do our interior elliptic
estimates for $u_j$. Now we write our equations (3.19) as follows:
 $$
\leqalignno{\text{(a)} \ \ & (D_j +E_j)u_j = u_j\widetilde H_j   &(3.20)\cr
           \text{(b)} \ \ & (D'_j +E_j)\bar u_j = \bar u_j\widetilde L_j.\cr}
$$

\n
Here the error operator $E_j = e_j^1 \f\pa{\pa z} + e^2_j\f\pa{\pa\bar z}$.
Let $a_j = u_j(p)$ and $\hat u_j = (a_j)^{-1}u_j$. Then $\hat u_j$
satisfies the same equation (3.20) as $u_j$ does but with $\hat u_j(p)=1$.
(note that $(m_{a_j^{-1}} \circ  \phi_j)^*J_j=\widetilde J_j$.  Let $U$ be a
small disc centered at $p$, which is compactly contained in another small
neighborhood $V$. Let $|\cdot|_{s,W}$ denote the $H^s$ norm on $W$. Let $D^*_j
\ (D^{'*}_j,$ resp.) denote the formal adjoint operator of $D_j + E_j\  (D'_j +
E_j$, resp.). It is easy to see that either $D^*_j \circ (D_j+E_j)$ or
$D^{'*}_j \circ (D'_j+E_j)$ is real positive self-adjoint, strictly and
uniformly elliptic in a neighborhood $\widetilde V$ of $p$ so that the
constants $\gamma$ and $\nu$ in (9.47) of [GT] are independent of $j$ for large
enough $j$. Choose small discs centered at $p, U, V_j, j=1,...,s, V$ such that
$U \subset V_1 \subset V_2 \subset\cdots\subset V_s \subset V \subset
\widetilde V$ where each smaller disc
is compactly contained in larger ones and $V$ is chosen so that we can apply
Theorem 9.20 of [GT]. By standard interior elliptic estimates, we compute (in
case $D=0$ at $p$, replace $D^*_j\circ (D_j+E_j)$ by $D^{'*}_j \circ
(D'_j+E_j)$ and $\hat u_j$ by $\overline{\hat u_j}$)
$$
\leqalignno{|\hat u_j|_{s+1,U} &\lesssim |D^*_j \circ (D_j+E_j)\hat
                            u_j|_{s-1,V_1} + |\hat u_j|_{0,V_1}   &(3.21)\cr
          &  \lesssim |D^*_j (\hat u_j\widetilde H_j)|_{s-1,V_1} + |\hat
             u_j|_{0,V_1}  \ \ \text{(by (3.20) (a))}\cr
          &  \lesssim |\hat u_j|_{s,V_1} \ \ \text{(by the interpolation
                 inequality)} \cr
          & \lesssim |\hat u_j|_{s-1,V_2}\ \ \text{(by the same argument as
                 above)}\cr
          &  ...\cr
          & \lesssim |\hat u_j|_{0,V} \cr}
$$

\n
where $A \lesssim B$ means $A \le k B$ for constant k independent of $\hat
u_j$.  On the other hand applying the Harnack estimates (Theorems 9.20, 9.22 in
[GT]) to the equation:
$$
   [D^*_j \circ (D_j+E_j)]\log |\hat u_j| = Re(D^*_j \widetilde H_j),
$$

\n
(noting that we apply theorems to $\log |\hat u_j|-\underset\,{\widetilde
V}\to\inf\log |\hat u_j| \ge 0$)  we obtain the estimate of the supremum norm
on $V$:
$$
|\hat u_j|_{L^\infty, V} \le C  \tag 3.22
$$

\n
where $C$ is a constant independent of large enough j. Combining (3.21) and
(3.22) we get
$$
|\hat u_j|_{s+1,U} \le  C_s,
$$

\n
so there exists a subsequence (still denoted $\{\hat u_j\})$ of $\{\hat u_j\}$
converging weakly in $H^{s+1}$ on $U$. By compactness, $\hat u_j$ converges in
any weaker norm, say, $L^2$ norm $|\cdot|_0$. By (3.21) with $U,V$ replaced by
$U',U (U'$ being a smaller disc centered at $p$, compactly contained in $U$)
resp., we learn that $\hat u_j$ is Cauchy in $H^{s+1}$ on $U'$. Therefore $\hat
u_j$ converges in $H^{s+1}$ on $U'$. Since $\hat N$ is compact, we can
pick up a finite number of such $(U', p)$ to cover $\hat N$. Let $p_i$'s denote
such points and $u_{j,U'_i}$ denote corresponding $u_j$ on $U'_i$. We adjust
$a_j$ and $\hat u_j$ to be $a_j = \max_i u_{j,U'_i}(p_i)$ and $\hat u_{j,U'_i}
= a_j^{-1} u_{j,U'_i}$. Thus $\hat u_{j,U'_i}(p_i) < 1$ so that our previous
argument still works  for all $i$. Now it is easy to pick a subsequence
of $m_{a_j^{-1}}\phi_j$, which converges in $H^{s+1}$ on each $U'_i$.
   \qquad   \qed

\bigskip
   Consider the action of $Bdiff_0^{s+1}$ or $\frak B^{s+1}$ on $\widetilde
P_{ic_s}$. First we describe the tangent space of the orbit passing through a
given element $\hat J$ in $\widetilde P_{ic}$. Push $\hat J$ down to a complex
structure $\hat c$ on $\hat N$. Take a local holomorphic coordinate $z$ of
$\hat N$ for $\hat c$. Take a local trivialization $(z,w)$ of $\hat L$ so that
$\pa_w$ and $Z_1 = \pa_z + b(z,\bar z)\bar w \pa_{\bar w}$ form a basis of the
type $(1,0)$ tangent vectors with respect to $\hat J$.
(note that $\hat J(\pa_z)=i\pa_z$ mod $\pa_w$ and $\pa_{\bar w}$,  and $Z_1$ is
$\Bbb C^*$-invariant) Want to find another trivialization $\widetilde z=z,
\widetilde w =\lambda(z,\bar z)w$  so that $\pa_w = \lambda \pa_{\widetilde w}$
and $Z_1 = \pa_{\widetilde z}$\  (mod$\ \pa_w)$. The chain rule tells us that
$$
   \pa_z = \pa_{\widetilde z} + \f{\pa\bar\lambda}{\pa
z}(\bar\lambda)^{-1} \bar w\pa_{\bar w}\ \ \text{(mod}\, \pa_w),
 $$

\n
so $\lambda$ has to satisfy the following $\bar \pa$-equation:
$$
\f{\pa\log\lambda}{\pa\bar z}= -\bar b.
$$

\n
But it is easy to solve the above equation locally. ($\lambda$ is in $H^{s+1}$
if $b$ is in $H^s$) Therefore we have a trivialization ($\widetilde
z,\widetilde w$) of $\hat L$, holomorphic with respect to $\hat J$, i.e.
$\{\pa_{\widetilde z}, \pa_{\widetilde w}\}$ forms a basis of type $(1,0)$
tangent vectors relative to $\hat J$. Now use $(z,w)$ instead of ($\widetilde
z,\widetilde w$) to denote a trivialization of $\hat L$, holomorphic with
respect  to $\hat J$, so $\hat J=i(dz \otimes \pa_z + dw \otimes \pa_w) $+
conjugate. Let $\phi_t$ be a family of $H^{s+1}$ bundle automorphisms of $\hat
L$. Recall that we write the infinitesimal bundle automorphism $V = \f
d{dt}|_{t=0} \phi_t = v^1 \pa_z + v^* w\pa_w$ + conjugate. (cf.(3.11)) Compute
$$
\f d{dt}|_{t=0}\phi_t^*\hat J = L_V \hat J
          = 2i\pa_{\bar z}v^1d\bar z \otimes \pa_z +
            2iw\pa_{\bar z}v^*d\bar z\otimes\pa_w + \text{conjugate}.
$$

\n
Recall that $\widetilde{\frak V}^s$ denote the Hilbert space of all
infinitesimal bundle automorphisms with bounded $H^s$-norm. Define the first
order operator $P: \widetilde{\frak V}^{s+1}\to T_{\hat J}\widetilde P_{ic_s}=
\frak E^s_{\hat J}$ by
$$
P(V) = L_V \hat J=2i v^1{}_{,\bar 1} d\bar z \otimes \pa_z + 2i v^*{}_{,\bar 1}
w d\bar z \otimes \pa_w +\text{conjugate}   \tag 3.23
 $$

\n
in trivialization $(z,w)$, where $v^1_{,\bar 1}= \pa_{\bar z}v^1, v^*_{,\bar
1}= \pa_{\bar z}v^*$. We want to describe the $L^2$ orthogonal subspace of
Range(P) in $\frak E^s_{\hat J}$,  which is supposed to be the kernel
Ker$P^*$ of the adjoint operator $P^*$. Since $\hat L$ is not compact, we
need to define a suitable inner product on $\frak E^s_{\hat J}$ over $\hat N$.
First observe from (3.7.1) that $E_1{}^{\bar 1}$ behaves just as a tensor
on $\hat N$ (under the special coordinate change (3.8)) while $E_1{}^*$ does
not by (3.7.2). We can adjust $E_1{}^*$ to get a tensor by the aid of
connection.  Let $\| \ \ \|$ be a hermitian metric on $\hat L$. Let $s(z)$
denote the local holomorphic section of $(\hat L,\hat J)$ given by $z\to(z,1)$
locally. Let $\nu =\|s(z)\|^2$. The canonical connection associated to $\|\ \
\|$ is given by
$$
\nu^{-1} \pa\nu = \Gamma(z,\bar z)dz
$$

\n
with $\Gamma = \pa_z \log\nu$. The transformation law according to (3.8)
goes as follows:
$$
\Gamma=\widetilde\Gamma h' + g'g^{-1}.
$$

\n
(noting that $g\widetilde s = s$) Define
$$
E_1 = E_1{}^* + E_1{}^{\bar 1}\bar\Gamma.   \tag 3.24
$$

\n
It is easy to see that $E_1 = \widetilde E_1 h'$, obeying the correct
transformation law as a tensor.  Let $g=g_{1\bar 1}dzd\bar z$ be the unique
hyperbolic metric on $\hat N$ associated to $\hat c$. We use $g^{1\bar 1}=
(g_{1\bar 1})^{-1}$ or $g_{1\bar 1}$ to raise or lower indices. Also denote
the volume form of g by dvol$_g$. Now we can define an inner product on
$\frak E^s_{\hat J}$:
$$
\langle E, F\rangle = \int_{\hat N}\{E_1{}^{\bar 1}F_{\bar 1}{}^1 + g^{1\bar
1}E_1F_{\bar 1}\} dvol_g.     \tag 3.25
 $$

\n
Here we have used the expression (3.6) for $E, F$ and (3.24) for $E_1$,
$F_{\bar 1}=\overline{(F_1)}$ . Take $E=P(V)$. Comparing (3.6) and (3.23) gives
$$
E_1{}^{\bar 1} = -2iv^{\bar 1}{}_{,1}, E_1{}^* = -2i\bar{v^*}{}_{,1}.\tag 3.26
$$

\n
Here $v^{\bar 1}=\overline{(v^1)}$  and $u_{,1} = \pa_zu$.  Define
$$
v = v^* + v^1\Gamma.    \tag 3.27
$$

\n
Easy to check that $v$ is independent of the choice of holomorphic
trivializations. Hence $v$ defines a global function on $\hat N$. Recall that
$c_1(=-m)$ denotes the first Chern number of $\hat L$. Let $\mu =
(1/4)c_1(\text{genus}(\hat N)-1)$.  For a special choice of $\|\ \ \|$
according to Proposition 2.1 relative to $\hat J$, we compute
$$
E_1 = -2i[\bar v_{,1} + \mu v_1]    \tag 3.28
$$

\n
where $v_1 = v^{\bar 1}g_{1\bar 1}$ and we have used $-g^{1\bar 1}\bar\Gamma_{
,1} = \mu$.  Substituting (3.26),(3.28) in (3.25) and using integration by
parts gives
$$
<P(V), F> = 2i \int_{\hat N}\{v^1(-F_1{}^{\bar 1}{}_{,\bar 1} +
\mu F_1)-vF_{1,\bar 1}g^{1\bar 1}\}dvol_g + \text{conjugate.}   \tag 3.29
$$

\n
(note that since $g$ is Kahler, the usual derivative of $v^1$ along $\bar z$
-direction coincides with its covariant derivative. Hereafter for a tensor $T$
on $\hat N, T_{,1} \ (T_{,1\bar 1}$ and so on, resp.) means the covariant
derivative of $T$ in the $z$-direction ($z\bar z$-direction and so on,
resp.))

The above formula suggests a suitable inner product on $\widetilde{\frak V}^s$
for our purpose. Namely, we define
$$
<V, U> = \int_{\hat N}[g_{1\bar 1}v^1 u^{\bar 1} + v\bar u] dvol_g +
\text{conjugate}    \tag 3.30
$$

\n
for $V= 2Re[v^1 \pa_z + v^* w \pa_w],\ U=2Re[u^1\pa_z + u^* w \pa_w]$ locally
and $v,u$ being global functions defined by (3.27). Define the adjoint operator
$P^*$ of $P: \frak E^s_{\hat J}\to\widetilde{\frak V}^{s-1}$ so that
$$
\langle P(V), F\rangle\ = \ \langle V, P^*(F)\rangle.
$$

\n
Then it follows from (3.29),(3.30) and (3.27) that locally
$$
P^*(F)= 2i(F_{\bar 1}{}^1{}_{,1} - \mu F_{\bar 1})g^{1\bar 1}\pa_z +
   2i[F_{\bar 1}{}^*{}_{,1}+F_{\bar 1}{}^*\Gamma+F_{\bar
1}{}^1(\Gamma_{,1}+\Gamma^2)]g^{1\bar 1}w\pa_w + \text{conjugate.}
$$

\n
If we represent $V$ by the pair $(v^1, v)$ and $E$ by the pair $(E_{\bar
1}{}^1, E_{\bar 1}$).  Then we can write $P(V)$ and $P^*(F)$ as follows:
$$
\leqalignno{& P(V) = 2i(v^1{}_{,\bar 1}, v_{,\bar 1}+\mu v_{\bar 1}) &(3.31)\cr
     & P^*(F) = 2i(g^{1\bar 1}(F_{\bar 1}{}^1{}_{,1}-\mu F_{\bar 1}), g^{1\bar
                1}F_{\bar 1,1}).  \cr}
$$

\n
Let $\Delta_{\hat J} = P^*P$. By (3.31), we compute
$$
\Delta_{\hat J}(V)= -4(g^{1\bar 1}(v^1{}_{,\bar 11}-\mu(v_{,\bar 1}+v_{\bar
1}\mu)),  g^{1\bar 1}(v_{,\bar 11} + \mu v_{\bar 1,1})).
$$

\n
The leading term of $\Delta_{\hat J}(V)$ is $-4(\Delta_gv^1, \Delta_gv)$
where the Laplacian $\Delta_g = g^{1\bar 1}\pa^2/\pa zd\bar z$. Thus
$\Delta_{\hat J}$ is a second order self-adjoint elliptic operator defined on
$\widetilde{\frak V}^s$.

\bigskip
\n
Lemma 3.6. Suppose $\hat J \in \widetilde P_{ic_{s+1}}$. Then there is an
$L^2$-orthogonal splitting
$$
\frak E^s_{\hat J} = Ker P^* + P(\widetilde{\frak V}^{s+1}).
$$

\n
Proof. Given $E$ in $\frak E^s_{\hat J}$. It is easy to see that $P^*(E)$ is
orthogonal to the kernel $Ker\Delta_{\hat J}$ of $\Delta_{\hat J}$ since
$Ker\Delta_{\hat J}=Ker P.$  Therefore by the standard elliptic theory we can
solve the equation  $\Delta_{\hat J}(V) = P^*(E)$ for $V$ in $H^{s+1}$. Now set
$E_0 = E - P(V)$. It is obvious that $E_0$ is in $Ker P^*$.
 \qquad \qed

\bigskip

   We remark that elements in $Ker P^*$ are all smooth by the elliptic
regularity. (note that $g_{1\bar 1}$ has the same regularity as $\hat J$ does
[Tr])

Moreover the dimension of $Ker P^*$ is finite. We compute it as follows. First
an element $F$ in $Ker P^*$ satisfies a system of linear equations:
$$
\leqalignno{&\ \ (a)\  F_1{}^{\bar 1}{}_{,\bar 1} - \mu F_1 = 0  &(3.32)\cr
            &\ \ (b)\  F_{1,\bar 1}= 0\cr}
$$

\n
by (3.31). Solutions $F_1$ for (3.32)(b) consist of all holomorphic (1,0)-forms
$F_1 dz$ on $\hat N$, denoted $H^{1,0}$. Let $Q(\hat N)$ denote the space
of holomorphic quadratic differentials on $\hat N$. By (3.32), the projection
map  from $Ker P^*$ onto $H^{1,0}$ has the kernel
equal to $Q(\hat N).$  From the basic linear algebra we learn that
$$
  \dim Ker P^* = \dim Q(\hat N) + \dim H^{1,0}.
$$

\n
On the other hand, $Q(\hat N)$ is known to describe the infinitesimal
Teichmuller space whose dimension is $6g-6$ by the Riemann-Roch theorem
(e.g. [Tr]) while dim $H^{1,0}$ is the same as that of the so-called
Picard variety in the Riemann surface theory, which is known to be $2g$.
Therefore
$$
\dim Ker P^* = 6g-6 + 2g = 8g - 6.   \tag 3.33
$$

\bigskip
\n
Lemma 3.7. Given $\hat J$ in $\widetilde P_{ic}$, there exists a local smooth
submanifold ${\frak S}$ of $\widetilde P_{ic_s}$ of dimension $8g-6$ passing
through $\hat J$ with the tangent space equal
to $Ker P^*$ at $\hat J$. Moreover, $\frak S$ consists of only smooth elements.

\bigskip
\n
Proof. Consider the map $\Phi_{\hat J}: Ker P^*\to \widetilde P_{ic_s}$. (see
(3.9)) It is easy to see that $\Phi_{\hat J}$ is smooth and its functional
derivative at 0 is the inclusion map from $Ker P^*$ into $\frak E^s_{\hat J}$,
which is surely injective and splits by Lemma 3.6.  Therefore $\Phi_{\hat J}$
is a smooth immersion at 0. That is to say, there exists a small neighborhood
$U$ of 0 such that $\frak S =\Phi_{\hat J}(U)$ is a smooth submanifold of
$\widetilde P_{ic_s}$ with dimension $8g-6$ by (3.33). Note that $\hat J$ is
smooth and elements in $Ker P^*$ are smooth as remarked previously. Thus $\frak
S$ consists of only smooth elements.
        \qquad\qed

\bigskip
   Now let $\Xi : \frak B^{s+1} \times \frak S\to \widetilde P_{ic_s}$ denote
the action of $\frak B^{s+1}$ on $\frak S$ by the pullback. Observe that $\Xi$
is smooth and
$$
D{\Xi}(id, \hat J): \frak V^{s+1} \times Ker P^*\to{\frak E}^s_{\hat J}
$$

\n
is given by $D{\Xi}(id, \hat J)([X], E) = L_X\hat J + E = P(X) + E$. If
$L_X\hat J=0$, then $X$ is an infinitesimal bundle automorphism fixing $\hat
J$. Thus $[X]=0$ by Proposition 3.4 and hence $D\Xi(id, \hat J)$ is a
continuous linear isomorphism by further using Lemma 3.6 and noting that
$\Delta_{\hat J}$ is elliptic.  Therefore $\Xi$ is a local diffeomorphism by
the inverse function theorem on Banach spaces. We have shown the existence of
"local slices":

\bigskip
\n
Proposition 3.8. There exist neighborhoods $W$ of $\hat J$ in $\widetilde
P_{ic_s}, U$ of id in $\frak B^{s+1}$ and $V$ of $\hat J$ in $\frak S$ such
that $\Xi: U \times V\to W$ is a diffeomorphism.

\bigskip
Now using freeness and properness of our $\frak B^{s+1}$ action
(Propositions 3.4, 3.5) plus the existence of "local slices" (Proposition 3.8),
we can equip our quotient space $\widetilde P_{ic}/{\frak B}$ with smooth
manifold structure by a standard argument. (e.g. section 2.4 in [Tr]) Recall
that we denote $\widetilde P_{ic}/{\frak B}$ by $P_{ic}^t.$

\bigskip
\n
Theorem 3.9. $P_{ic}^t$ is a smooth manifold of dimension $8g-6.$

\bigskip
\n
Proof. First we show the existence of "slices": that is to say, if we take
the slice $\frak S$ to be sufficiently small, then each orbit of $\frak B$
passing through $\frak S$ intersects $\frak S$ at exactly one point, i.e.
$\phi^*J$ in $\frak S$ with $J$ in $\frak S$ implies $\phi=id$. Suppose this
is not true. Then there are sequences $\phi_j$ in $\frak B$ and $J_j$ in
$\frak S$ such that $J_j$ and $\phi_j^*J_j$ converge to $\hat J$ in $H^s$
while all $\phi_j$'s keep ontside some fixed $H^{s+1}$ neighborhood of
$id$ in $\frak B$ in view of Proposition 3.8. (we equip $\widetilde P_{ic},
\frak B$ with the $H^s, H^{s+1}$ topologies, resp.) By Proposition 3.5
(properness) there exists a subsequence of $\phi_j$, which converges to $\phi$
in $H^{s+1}$.  It follows that $\phi^*\hat J =\hat J$ and then $\phi=id$ by
Proposition 3.4 (freeness), contrary to $\phi_j$'s sitting outside some
neighborhood of $id$. Thus we can take the slices as coordinate charts (instead
of their tangent spaces). It is easy to see by Proposition 3.8 that the
transition function is smooth.
\qquad                         \qed

\bigskip
\n
Proof of Theorem C: We will introduce a natural complex structure on
$P_{ic}^t$. First there is a canonical way to define an almost complex
structure $\Theta$ on $\widetilde P_{ic}$: for $J$ in $\widetilde P_{ic}, E$ in
${\frak E}_J$,
$$
    \Theta_J(E) = J \circ E.
$$

\n
It is easy to verify that $J \circ E$ is still sitting in $\frak E_J$.  Let
$\pi: \widetilde P_{ic}\to P_{ic}^t$ be the natural projection. From our
previous argument $({\pi}, {\frak B}, \widetilde P_{ic}, P_{ic}^t)$ is a
(weak) principal ${\frak B}$-bundle in the sense of [Tr], p.54. (note that the
right action of $\frak B$ on $\widetilde P_{ic}$ is given by pulling back)
It is straightforward that $\Theta$ is $\frak B$-invariant (cf. p.88 in [Tr]),
and $\Theta$ maps "vertical" vectors to "vertical" vectors: (see p.86 in [Tr]
for the definition) since each J in $\widetilde P_{ic}$ is integrable by
Proposition 3.1, the associated Nijenhuis tensor
vanishes. It follows that ${\Theta}_J(L_XJ)=J L_XJ=L_{JX}J$ (cf. p.88 in
[Tr]), so $\Theta$ makes ($\pi, \frak B, \widetilde P_{ic}, P_{ic}^t)$
into an almost complex principal $\frak B$-bundle. (see Definition 4.1.4 on
p.86 in [Tr])  Next we note that the Lie bracket of two vector fields on
$\widetilde P_{ic}$ can be defined as in [Tr], p.85: instead of using
projections, we view $DY(J)X(J) =d/dt|_{t=0} Y(J(t))$ with $J(0)=J,
J'(0)=X(J);$ verify $DY(J)X(J)-DX(J)Y(J)$ is in
$\frak E_J$ for $X(J), Y(J)$ in $\frak E_J$ by observing that an element $E$ in
$\frak E_J$ can be described by the following conditions:
$$
\align
  &   E \circ J + J \circ E = 0\\
  &   m_\rho^* E = E\\
  &   E(v)=0\ \ \text{for}\ \ v \ \ \text{tangent to fibres of}\ \ \hat L.
\endalign
$$

Now we can define the Nijenhuis tensor $N({\Theta})$ of $\Theta$ on $\widetilde
P_{ic}$ as usual. Then a direct computation as shown in [Tr], p.88 yields
$N(\Theta)=0.$  By Theorem 4.1.2 in [Tr], the almost complex structure
$J_{pic}$ on $P_{ic}^t$  induced from $\Theta$ on $\widetilde P_{ic}$
has the vanishing Nijenhuis tensor. Since
$P_{ic}^t$ is a finite dimensional manifold, $J_{pic}$ is integrable,
i.e. there exists a complex structure on $P_{ic}^t$ whose associated almost
complex structure is $J_{pic}$ by the Newlander-Nirenberg theorem.
  \qquad\qed

\bigskip
\n
4. A supporting manifold of $\frak S^t$ and properness of the contact
action.

   Recall (cf. section 1) that $\frak S^t$ is the quotient space of $\widetilde
{\frak S}$ modulo $C^0_{\hat H}$. Here $\widetilde{\frak S}$ denotes the space
of  all smooth spherical CR manifolds $(\hat S, \hat H, J)$ with $J$ oriented
and compatible with $\hat H$ and $C^0_{\hat H}$ denotes  the identity component
of the orientation-preserving smooth contact
diffeomorphism group $C_{\hat H}$ relative to $\hat H$. In this section we will
parametrize a local "supporting" space of $\frak S^t$ and show the
properness of the $C^0_{\hat H}$ action.

   We will work with the aid of anisotropic Folland-Stein spaces. For
$F$ a vector bundle over a closed contact manifold $(M,H)$ and $k$ a
nonnegative integer, let $S_k(F)$ denote the $L^2$  Folland-Stein space of
sections of $F$. ([FS], p.241 in [CL1]) If the bundle is clear from the
context, we simply use the notation $S_k$ instead of  $S_k(F),$ and a norm on
$S_k$ is denoted by  $|\cdot|_k$. Let $\widetilde{\frak S}^k$ denote the
completion of $\widetilde{\frak S}$ under the norm $|\cdot|_k$ for a fixed
smooth background contact  manifold
$(\hat S, \hat H)$. Let $\frak J^k (\frak J$,resp.) denote the space of all
oriented compatible $S_k\ (C^\infty$,resp.) CR structures on ($\hat S, \hat
H$) or a general contact manifold $(M,H)$ depending on the context. (note that
these CR  structures are sections of the endomorphism bundle End ($\hat H$))

\bigskip
\n
Lemma 4.1. Suppose dim $M=3$. For $k \ge 6$, (a) $S_k$ is an algebra;
(b) Let $f$ be a smooth function on nonnegative real numbers. Then $f\circ h$
is  still in $S_k$ for nonnegative $S_k$ function $h$.

\bigskip
\n
Proof. (a) is known. (e.g. [BD]) (b) is probably also known. We prove it
by induction on $k$. Computing the derivative of $f \circ h$ in some contact
direction give the derivative of $f$ composed with $h$ times the derivative of
$h$ in that direction, so induction hypothesis on $k-1$ plus (a) implies the
derivative of $f \circ h$ is in $S_{k-1}$. Hence $f \circ h$ is in $S_k$, so to
complete the proof we have to check the starting case $k=6$. But it is
straightforward by observing that $S_6$ is contained in $S_1^{12}$ or $S_2^8$
and $S_3$ is contained in $C^0$, etc.. (see e.g. Theorem 4.17, Corollary 5.16
in [Fo])                \qquad\qed

\bigskip
\n
Lemma 4.2. $\frak J^k$ is a Hilbert manifold for large $k$, say, $k\ge 6$.

\bigskip
\n
Proof. In [CL1], we parametrize $\frak J^k$ for $k=\infty$, i.e. in the smooth
category by a map $\Phi_J$ given by
$$
  \Phi_J(E) = E_0 J + E,\ \ \   E_0 = (1 + (1/2)Tr(E^2))^{1/2}.
$$

\n
(p.228, Lemma 2.3 in [CL1]) Suppose $E$ is in $S_k$.  (more precisely
$S_k({\frak E}_J))$ By Lemma 4.1(a) $h = E_1{}^{\bar 1}E_{\bar 1}{}^1$ is in
$S_k$. Take $f(x) = (1 + x)^{1/2}$. By Lemma 4.1(b) $E_0$ is
in $S_k$. Therefore $\Phi_J$ preserves $S_k$ spaces, so does $\pi_J$. Thus
we can still use $\Phi_J$ or $\pi_J$ to parametrize $\frak J^k$ modelled
on $S_k(\frak E_J)$.\qquad\qed

\bigskip
   Hereafter throughout this paper we will assume that $k \ge 6$ unless
specified otherwise. We know that a CR structure $J$ being spherical is
characterized by the vanishing of the Cartan tensor $Q_J$. (p.227 in [CL1]) The
linearization $DQ_J$ is subelliptic when restricted to Ker $B_J$. (in view of
Lemma 3.3 and Proposition 3.1 in [CL1])
When working in the Folland-Stein category, it is enough to still require
the reference CR structure to be smooth for our purpose.

\bigskip
\n
Lemma 4.3. For a smooth spherical $\hat J$ in $\frak J^k$ (and an auxiliary
smooth contact form), we have the following $L^2$-orthogonal decomposition:
$$
S_k(\frak E_{\hat J}) = Ker_k DQ_{\hat J} + DQ_{\hat J}(S_{k+4}(\frak E_{\hat
J}))    \tag 4.1
$$

\n
where $Ker_k$ means elements in the kernel and also in $S_k.$

\bigskip
\n
Proof. Differentiating the Bianchi identity  $B_JQ_J=0$ in Proposition 3.1 of
[CL1] at $\hat J$ in the direction E implies $DQ_{\hat J}(E)$ belongs to the
kernel of $B_{\hat J}$.  (note that $Q_{\hat J}=0$)  On the other hand, for $E$
in $Ker_kB_{\hat J}$, we have
$$
   DQ_{\hat J}(E) = -(1/24)L_a^* L_a(E) + \ \text{terms of lower weight}
$$

\n
with $a = 4 + i\sqrt 3$ according to Lemmas 3.3, 3.2 in [CL1]. For $a$ not
an odd integer, $L_a$ is a subelliptic operator of weight 2, i.e. satisfies
the estimate (4.2) in [CL1], so restricted to $Ker_kB_{\hat J}, DQ_{\hat J}$ is
a subelliptic operator of weight 4 according to the above formula, i.e. earns
four derivatives in contact directions and we have the $L^2$-orthogonal
decomposition for $DQ_{\hat J}$:
$$
 Ker_k B_{\hat J} = Ker DQ_{\hat J} + DQ_{\hat J}(Ker_{k+4} B_{\hat J}).
              \tag 4.2
$$

\n
Here $Ker DQ_{\hat J}$ consists of smooth elements since $DQ_{\hat J}$ is
subelliptic, hence hypoelliptic when restricted to $Ker_k B_{\hat J}$. We also
have the $S_k$ version of Proposition 2.4 in [CL2]:(note that notation $D_J$ in
[CL2]= $B'_J$ in [CL1])
$$
 S_k(\frak E_{\hat J}) = Ker_k B_{\hat J} + B'_{\hat J}(S_{k+2})    \tag 4.3
$$

\n
basically because $\Delta_{\hat J}= B_{\hat J}B'_{\hat J}$ is a subelliptic
operator of weight 4 by Lemma 2.1 in [CL2]. Since each element in the range of
$B'_{\hat J}$ is an infinitesimal contact orbit at $\hat J$ and $Q_J$ equals 0
for $J$ in the contact orbit of $\hat J, DQ_{\hat J}$ vanishes on $B'_{\hat
J}(S_{k+2})$. Therefore we can combine (4.3) and (4.2) to get (4.1).
\qquad\qed

   Take a smooth $\hat J$ in $\widetilde{\frak S}$ and choose an auxiliary
smooth contact form $\hat\theta$. There is a local slice $\frak S$ of $\frak J$
passing through $\hat J$ by Theorem A of [CL2], defined by $\Phi_{\hat J}(Ker
B_{\hat J}$) (note $B_{\hat J}=D^*_{\hat J})$ for elements in Ker $B_{\hat J}$
with small $|\cdot|_{5,\infty}$ norm. By the Sobolev lemma for our anisotropic
spaces (e.g., (4.17), (5.15) in [F]), we have $S_k \subset S_{8-4/q} \subset
S_6^q \subset {\Gamma}_{6-4/q} \subset S_5^\infty$ for $k \ge 8$.
Thus taking elements of small $S_k$ norm, $k\ge 8$, in $Ker\ B_{\hat J}$ and
then sending them to ${\frak J}^k$ through $\Phi_{\hat J}$, we obtain an $S_k$
slice $\frak S_{(k)}$ passing through ${\hat J}.$  Consider the map $\frak Q:
\frak S_{(k)}\to DQ_{\hat J}(S_k(\frak E_{\hat J})$, defined by
$$
   \frak Q(J) = \pi(Q_J).
$$

\n
Here $Q_J$ is the Cartan tensor of $(J, \hat\theta), \pi$ is the composition
of the orthogonal projection $\pi_{\hat J}: S_{k-4}(End(\hat H))\to
S_{k-4}(\frak E_{\hat J})$ (p.228 in [CL1]) and the projection: $S_{k-4}(\frak
E_{\hat J})\to DQ_{\hat J}(S_k(\frak E_{\hat J}))$ according to (4.1).

\bigskip
\n
Proposition 4.4. $\frak Q^{-1}(0)$ is a smooth finite dimensional submanifold
of $\frak S_{(k)}$ for $k\ge 10$ near a smooth $\hat J$.

\bigskip
\n
Proof. It is easy to see that $\pi$
is smooth and since $Q_J$ is of type 4 (pp.249-250 in [CL1]), the map: $J\in
S_k\to Q_J\in S_{k-4}$ for $k \ge 10$ is smooth. (note that $S_k$ forms an
algebra for $k \ge 6$ by Lemma 4.1) Therefore $\frak Q$ is smooth. We compute
$$
\leqalignno{D\frak Q(\hat J)(E) &= D\pi(0)DQ_{\hat J}(E)   &(4.4)\cr
                   & = \pi(DQ_{\hat J}(E)) = DQ_{\hat J}(E) \cr}
$$

\n
for $E$ in $S_k(\frak E_{\hat J})$.  From (4.4) it is clear that $D\frak Q(\hat
J)$ is surjective. Furthermore, the kernel of $D\frak Q(\hat J)$ is the same
as the kernel of $DQ_{\hat J}$, which splits according to (4.1). Thus by the
inverse function theorem $\frak Q$ is a submersion at $\hat J$ (Proposition 2
on p.27 in [La]), so $\frak Q^{-1}(0)$ has a smooth submanifold structure
near $\hat J$. Moreover, finite-dimensionality follows from subellipticity of
$DQ_{\hat J}$ restricted to $Ker B_{\hat J}$.
      \qquad           \qed

\bigskip
   We will use $\frak Q^{-1}(0)$ as a "supporting" background
manifold to prove $\frak S^t_0$ (an open connected subspace of
$\frak S^t$; see section 5) is a smooth manifold. First we will show the
properness of the contact action in the negative pseudohermitian curvature
case. Let $(M,H)$ be a smooth, closed, oriented, contact 3-manifold.
Let $S_k(M,\Bbb R)$ denote the space of all real-valued $S_k$ functions on
$(M,H)$.

\bigskip
\n
Lemma 4.5. The pseudohermitian curvature $R_{J,\theta}$ belongs to
$S_{k-2}(M,\Bbb R)$ for $J, \theta$ in $S_k$ with $k \ge 8$.

\bigskip
\n
Proof. Take a smooth contact form $\hat\theta$. Write $\theta = u^2 \hat\theta$
for positive $u$ in $S_k(M,\Bbb R)$. The transformation law reads
$$
 4\Delta_bu + R_{J,\hat\theta}u - R_{J,\theta}u^3 = 0.  \tag 4.5
$$

\n
([JL]) Here the negative sublaplacian $\Delta_b$ is defined with respect to
($J, \hat\theta$). Suppose $J$ is in $S_k$. Then it is easy to see that
$\Delta_bu$ is in $S_{k-2}$ in view of Lemma 4.1 if we write $J$ with respect
to a smooth $\hat J$ as in p.249 of [CL1] and apply formulas on pp. 249-250 of
[CL1] to express $\Delta_bu$. Moreover $R_{J,\hat\theta}$ is in $S_{k-2}$ since
$R_{J,\hat \theta}$ is of type 2 as shown in the following transformation
formula:
$$
\leqalignno{R_{J,\hat\theta} = & R_{\hat J,\hat\theta} + \f 12 i(v_{\bar
1}{}^1{}_{,0}v_1{}^{\bar 1}-v_1{}^{\bar 1}{}_{,0}v_{\bar 1}{}^1) &(4.6)\cr
   & - v_0(v_{0,1\bar 1}+v_{0,\bar 11} + v_{11,\bar 1\bar 1} + v_{\bar
          1\bar 1,11})\cr
   & - v_{\bar 1\bar 1}(v_{0,11} + v_{11,\bar 11}) -v_{11}(v_{0,\bar 1\bar
      1} + v_{\bar 1\bar 1,1\bar 1})-2|v_{0,1} + v_{11,\bar 1}|^2  \cr}
$$

\n
(see pp.249-250 in [CL1] where we did not give the above formula precisely)
Now from (4.5) $R_{J,\theta}$ is therefore in $S_{k-2}$ in view of Lemma 4.1.
         \qquad\qed

\bigskip
\n
Lemma 4.6. Let $(M,H)$ be a smooth, closed, oriented, contact 3-manifold.
Suppose the pseudohermitian curvature $R_{\hat J,\hat\theta} =-1$ for some
smooth $(\hat J,\hat\theta)$ on $(M,H)$. Then for any $J$ in $\frak J^k,\ k\ge
8$, there exists a uniquely determined $S_k$ contact form $\theta$ with
$R_{J,\theta} = -1.$

\bigskip
\n
Proof. Consider the map $\frak R: \frak J^k \times \{S_k$ contact forms$\}\ \to
S_{k-2}(M,\Bbb R)$ defined by
$$
   {\frak R}(J,\theta) = R_{J,\theta}.
$$

\n
(well defined by Lemma 4.5) The map $\frak R$ is smooth in view of (4.5) and
(4.6). Differentiating $\frak R$ at $(\hat J,\theta)$ in the direction
$(J',\theta')=(2E,2h\hat\theta)$ gives
$$
D\frak R(\hat J,\hat\theta)(2E,2h\hat\theta)
= i(E_{11,\bar 1\bar 1} - E_{\bar 1\bar 1,11}) - (A_{11} E_{\bar
          1\bar 1} + A_{\bar 1\bar 1} E_{11})  + 4 {\Delta}_bh - 2 R_{\hat
          J,\hat\theta}h   \tag 4.7
$$

\n
according to (2.20) in [CL1] and (5.15) in [Lee]. Since $R_{\hat
J,\hat\theta} = -1, 4\Delta_b - 2R_{\hat J,\hat\theta} = 4\Delta_b + 2Id$. is
invertible. It follows that $D\frak R(\hat J,\hat\theta)$ is surjective.
Moreover it is easy to see that the kernel $(D\frak R(\hat
J,\hat\theta))^{-1}(0)$ and the space $\{(0,2h\hat\theta)\}$ span the tangent
space of the domain at $(\hat J,\hat\theta)$ and have only (0,0) in their
intersection. That is to say, $(D\frak R(\hat J,\hat\theta))^{-1}(0)$ splits.
Therefore $\frak R$ is a submersion at $(\hat J,\hat\theta)$. (Prop.2 on p.27
in [La])  Thus $\frak R^{-1}(-1)$ has a submanifold structure near ($\hat
J,\hat\theta)$ and it projects onto a neighborhood of $\hat J$ in $\frak J^k.$

   On the other hand, suppose $R_{J_j,\theta_j} = -1$ for a sequence of
smooth $(J_j,\theta_j$). (note that $\theta_j$ is uniquely
determined by $J_j$ by Prop. 2.2)  If $J_j$ tends to $J$ in $S_k$, we claim
that $\theta_j$ tends to $\theta$ in $S_k$ too so that $R_{J,\theta}
= -1$. Write $\theta_j = u_j^2 \hat\theta$ for positive $u_j$.
Then $u_j$ satisfies the equation (4.5) with $(J,\theta)$ replaced by
$(J_j,\theta_j).\ \ R_{J_j,\theta_j} = -1$ implies $R_{J_j,\hat\theta}$
must be negative by the maximum principle.  Moreover apply the maximum
principle to the equation (4.5) where $u_j$ is a maximum, hence (negative
sublaplacian) ${\Delta}_bu_j\ge 0$. Since both $R_{J_j,\hat\theta}$ and
$R_{J_j,\theta_j}$ are negative, we get the uniform $C^0$ estimate of
$u_j$:
$$
\max u_j \le (-R_{J_j,\hat\theta})_{\max}^{1/2} \le C   \tag 4.8
$$

\n
for a constant $C$ independent of $j$ in view of (4.6). Similarly applying
the maximum principle at the minimum of $u_j$, we obtain
$$
 0 \le c \le [(-R_{J_j,\hat\theta})_{\min}]^{\f 12} \le \min u_j   \tag 4.9
$$

\n
for a positive constant $c$ independent of $j$. Let $\Delta_b$ and
${\Delta}_{b(j)}$ denote the negative sublaplacians with respect to
($J,\hat\theta$) and $(J_j,\hat\theta)$ respectively. Using those formulas on
pp.249-250 in [CL1], we have the following estimate: given a small $\epsilon >
0,$
$$
 |\Delta_{b(j)}u -\Delta_bu|_{k-2} \le \epsilon |u|_k   \tag 4.10
$$

\n
for $j$ large and $u$ in $S_k$. For $J$ in $S_k$ the difference between
${\Delta}_b$ and the corresponding operator on the Heisenberg group is small
for a small region on $M$ in the sense of (4.10). By a standard argument
(absorbing the right side of (4.10) and using a partition of unity for compact
M), we still have the subelliptic estimate for $\Delta_b$:
$$
 |u_j|_k \le  C(|\Delta_bu_j|_{k-2} + |u_j|_0).   \tag 4.11
$$

\n
Write $\Delta_bu_j=\Delta_{b(j)}u_j + (\Delta_b-\Delta_{b(j)})u_j$
and substitute in (4.11). Using (4.10), absorbing the right side to the
left and applying the equation (4.5) to $\Delta_{b(j)}u_j$, we obtain
$$
 |u_j|_k \le C'.     \tag 4.12
$$

\n
Here $C'$ is a constant independent of $j$ and we have used (4.8) in estimating
$\Delta_{b(j)}u_j$ and dominating the $L^2$ norm of $u_j$. From (4.12) there
exists a subsequence, still denoted $u_j$, which weakly converges to $u$ in
$S_k$ but strongly converges to $u$ in $S_{k-1}$, say. Let $\theta=u^2
\hat\theta.$  Applying (4.11) to $u_j - u$ and using the interpolation
inequality to absorb $|u_j-u|_{k-2}$ to the left side, we get
$$
|u_j-u|_k \lesssim |R_{J_j,\hat\theta}-R_{J,\hat\theta}|_{k-2}+|u_j-u|_3.
          \tag 4.13
$$

\n
Here we have used an interpolation inequality (Cor.2.11 in [BD]) to estimate
$u_j^3 - u^3$. It follows by (4.13) that $u_j$ tends to $u$ in $S_k$, and it
is clear that $R_{J,\theta} = -1$ in view of (4.5),(4.6). We have proved
our claim.  Now consider the space $\frak J_{-1}$ of all smooth $J$ in $\frak
J$ such that $R_{J,\theta} = -1$ for some smooth $\theta$. (unique if
exists) The argument in our first paragraph shows that ${\frak J}_{-1}$ is open
in ${\frak J}$ (in $C^\infty$ topology). The argument (and our claim) above
shows in particular that $\frak J_{-1}$ is closed in $C^\infty$ topology.
Therefore $\frak J_{-1}= \frak J$. The lemma follows since $\frak J^k$ is the
completion of $\frak J$ under the norm $|\cdot|_k.$      \qquad\qed

\bigskip
   We remark that the similar idea of the above proof has been applied
to the case of a fixed CR structure in [CH]. We can now prove the properness
of contact diffeomorphisms acting on $\frak J$ in the case of negative
pseudohermitian scalar curvature. We can talk about $S_k$ contact
diffeomorphism on a contact manifold. (see Prop.2.18 in [BD])

\bigskip
\n
Proposition 4.7. Let the assumptions be as in Lemma 4.6. Let $\phi_j$ be a
sequence of contact diffeomorphisms in $S_{k+1}$ with $k \ge 12$. Suppose
$\phi_j^*J_j$ and $J_j$ converge in $S_k$ as $j$ goes to infinity for
$J_j$ in $\frak J^k$. Then there exists a subsequence of $\phi_j$ which
converges in $S_{k+1}$.

\bigskip
\n
Proof. From Lemma 4.6 we can associate a unique $S_k$ contact form
$\theta_j$ to $J_j$ so that $R_{J_j,\theta_j} = -1$. Let $g_j$ be the
adapted metric associated to $(J_j,\theta_j)$:([CH]) i.e. $g_j =
\theta^2_j + d\theta_j(\cdot,J_j(\cdot))$.
$g_j$ converges at least in $S_{k-2}$ since $\theta_j$ converges in $S_k$
as shown in the proof of Lemma 4.6. $\phi_j^*\theta_j$ is just
the unique contact form associated to $\widetilde J_j = \phi_j^*J_j$
satisfying the equation of pseudohermitian scalar curvature = -1. It follows
that $\phi_j^*g_j$ converges at least in $S_{k-2}.\  S_{k-2}$ is
contained in $H^{(k-2)/2}$ (the usual $L^2$ Sobolev space) with $(k-2)/2 >
4$. Therefore we can apply the result of Ebin and Palais (Theorem 2.3.1 in
[Tr]) to conclude the convergence of a subsequence (still denoted $\phi_j$) of
$\phi_j$ in $H^{k/2}$. We need to show the convergence actually is in
$S_{k+1}$. Take a smooth contact form $\hat\theta.$  There is a uniquely
determined smooth vector field $\hat T$ such that $\hat\theta(\hat T) = 1,
d\hat\theta(\hat T,\cdot) = 0. $ For $(J_j,\hat\theta)$ we can choose $S_k$
admissible coframe $\theta^1_j$. ([Lee])
(let $e_1$ be a smooth local section of the contact bundle $H$. Let $\omega^1,
\omega^2, \hat\theta$ be a local coframe dual to $e_1, Je_1, \hat T$. Then
$\theta^1$ is defined to be $\omega^1+ i\omega^2$ and if $J$ is in $S_k$, then
$\omega^1, \omega^2$, hence $\theta^1$ is in $S_k$) Write $\theta_j =
e^{2g_j} \hat\theta.\ \ g_j$ converges in $S_k$ since $\theta_j$ converges in
$S_k$. Also write $\phi_j^*\theta_j = e^{2f_j}\hat\theta.\ \ f_j$ converges in
$S_k$ since $\widetilde J_j$ converges in $S_k$ by the assumption. (same
reasoning as in the proof of Lemma 4.6) It is easy to see
$$
\phi_j^*\hat\theta = e^{-2g_j \circ \phi_j + 2f_j}\hat\theta   \tag 4.14
$$

\n
Let $h_j = -g_j \circ \phi_j + f_j$.  Let $\widetilde\theta^1_j$ be a local
$S_k$ admissible coframe with respect to ($\widetilde J_j,\hat\theta$). Then we
can adjust $\theta^1_j$ in $S_k$ by a modulus 1 factor (still denoted
$\theta^1_j$) so that
$$
\phi_j^*\theta^1_j = e^{h_j}\widetilde\theta^1_j \ \ \text{modulo}\ \hat\theta.
     \tag 4.15
$$

\n
((5.5) on p.421 in [Lee]) Now suppose $\phi_j$ converges in $S_l$ for $l \le
k.$  Then $h_j$ converges in $S_l$ too. (the composition map of an $S_l$
function and an $S_l$ contact diffeomorphism is still $S_l$ and the map is
jointly continuous for $l \ge 6$. A proof can be given by mimicking the one for
the usual $L^2$ Sobolev spaces. See pp.15-16 in [Eb]. Also see Prop. 2.13 in
[BD] for the precise estimate. Note that we start with $S_{k/2}$ with $k/2 \ge
6$  in which $\phi_j$ converges) It follows that $\phi_j^*\theta^1_j$ and
$\phi_j^*\theta^{\bar 1}_j$ converge in $S_l$ when applied to vectors tangent
to the contact bundle by (4.15). Let $\widetilde Z_1 = e_1 + i\widetilde Je_1$
where $\widetilde J$ is the limit of $\widetilde J_j$ in $S_k$.
Then $\phi_{j*}(\widetilde Z_1)$ and $\phi_{j*}(\widetilde Z_{\bar 1}$)
converge in $S_l$. Therefore $\phi_j$ converges in $S_{\ell+1}$. Thus by
induction we finally obtain that $\phi_j$ converges in $S_{k+1}$.
\qquad\qed

\bigskip
We remark that the properness of the contact action for a contact manifold
is generally not true. For instance, say, $\frak J$ contains a CR structure
with noncompact CR automorphism group. Now we can apply Proposition 4.7 to our
case $(M,H)=(\hat S,\hat H)$ on which there are canonical spherical ${\hat J}$
and contact form $\hat\theta$ such that $R_{\hat J,\hat\theta} = -1$. (and
$A_{\hat J,\hat\theta} = 0$. See section 2)  Let $C_J$ denote the group of CR
automorphisms relative to $J$ with the identity component $C^0_J$. We have a
$U(1)$ action on $\hat S$ given by fibre multiplications
by unit-length constants. Let $\widetilde{\frak S}^{U(1)}$ denote the space of
$U(1)$ invariant elements in $\widetilde{\frak S}$.  Let $C^{U(1),0}_{\hat H}$
denote the identity component of the group of $U(1)$ equivariant contact
diffeomorphisms in $C^0_{\hat H}$.

\bigskip
\n
Proposition 4.8. (a) For $J$ in $\widetilde{\frak S}^{U(1)}, C^0_J$
equals $U(1)= \{$fibre multiplications by unit-length constants$\}$ and is
contained in the center of $C^{U(1),0}_{\hat H}.$

\bigskip
(b) $C^{U(1),0}_{\hat H}/U(1)$ acts on $\widetilde{\frak S}^{U(1)}$
freely and properly.

\bigskip
\n
Proof. Any CR automorphism $\phi$ in $C^0_J$ relative to $J$ in $\widetilde
{\frak S}^{U(1)}$ is $U(1)$-equivariant by Proposition 3.14 in [Ep].
Therefore it can be pushed down to a biholomorphism on $\hat N$, which
must be the identity since genus $(\hat N) \ge 2$. On the other hand $\phi$
extends to  a holomorphic bundle automorphism of $\hat L$.  Therefore $\phi$
is just a fibre multiplication by a nonzero holomorphic function on $\hat N$,
which must be constant since $\hat N$ is closed (compact without boundary).
(cf. Proposition 3.4) The second conclusion of (a) follows by the definition of
$C^{U(1),0}_{\hat H}.$  Now (b) is clear by (a) and Proposition 4.7.
    \qquad\qed

\bigskip
   In the next section we will parametrize a certain open connected subspace
$\frak S^t_0$ of $\frak S^{t,U(1)}=\widetilde {\frak S}^{U(1)}/
C^{U(1),0}_{\hat H}$  as a smooth manifold and show that $\frak S^t_0$ is
diffeomorphic to $P_{ic_0}^t$, an open connected subspace of $P_{ic}^t.$

\bigskip
\n
5. The smooth manifold structure on $\frak S^t_0$: Proof of Theorem A
   and Corollary B

   Let $\widetilde P_{ic_0}$ be the connected component of $\widetilde P_{ic},$
containing ($\hat L,\hat J$). Define $P_{ic_0}^t$ to be the quotient
space of $\widetilde P_{ic_0}$ modulo the action of $\frak B$ or $Bdiff_0.\ \
P_{ic_0}^t$ is an open connected subset of $P_{ic}^t$. (actually
they are the same since $\widetilde P_{ic}$ is known to be connected.
But we do not pursue it here)

   Given an element ($\hat L,\widetilde J$) in $\widetilde P_{ic_0}$, there
associates a unique (up to a positive constant multiple) hermitian metric $\|
\ \ \|_{\widetilde J}$ on $\hat L$ according to Proposition 2.1. Define $\rho:
\hat L\backslash$the zero section $\to \Bbb R$ by $\rho(s)=\|s\|_{\hat
J}/\|s\|_{\widetilde J}$. Here $\hat J$ denotes the complex structure on $\hat
L$ (and also $\hat S$) associated to the fixed (or reference) holomorphic line
bundle ($\hat L,\hat N$) as before. It follows that $\rho(\lambda s)=\rho(s)$
for $\lambda$ in $\Bbb C\backslash\{0\}$, so $\rho$ can be pushed down to
define a function on $\hat N,$ still denoted $\rho$. Define $m_\rho:\hat
L\to\hat L$ by
$$
   m_\rho(s)=\rho(\hat\pi(s))s
$$

\n
where $\hat\pi:\hat L\to\hat N$ is the projection. Note that $m_\rho$ maps
$\hat S= \{s \in\hat L : \|s\|_{\hat J}=1\}$ onto $S_{\widetilde J}=\{s \in
\hat L: \|s\|_{\widetilde J}=1\}$. The contact bundle $\widetilde H$ defined by
subbundle  of $T\hat S$, invariant under the endomorphism $m_\rho^*\widetilde
J$ restricted to $T\hat S$, differs from $\hat H$ in general. We need
the $U(1)$-invariant version of Gray's theorem. Let $M$ be a closed (compact
without boundary) smooth manifold of dimension $2n+1$ with a smooth $U(1)$
action.  Suppose for each $\xi$ in $U(1)$, the action $A_\xi$ on $M$ is a
diffeomorphism.  Let $Diff^{U(1)}(M)$ denote the space of all
$U(1)$-equivariant diffeomorphisms.  Let $\frak B^{U(1)}$ denote the space of
all $U(1)$-invariant (smooth) contact bundles. It is clear that
$Diff^{U(1)}(M)$ acts on $\frak B^{U(1)}$ by pushing
forward. In Appendix A, we will show that both $Diff^{U(1)}(M)$ and
$\frak B^{U(1)}$ are smooth tame Frechet manifolds in the terminology of [Ha];
we will also show the following $U(1)$-invariant version of Gray's theorem.
(cf. Theorem 2.4.6 in [Ha])

\bigskip
\n
Theorem 5.1. Any contact bundle near a given one $H$ in $\frak B^{U(1)}$ is
conjugate to $H$ by a $U(1)$-equivariant diffeomorphism near the identity.
The identity component of $Diff^{U(1)}(M)$ acts transitively on each component
of $\frak B^{U(1)}$.

\bigskip
   Now apply Theorem 5.1 to our case: $M=\hat S$ with the $U(1)$ action given
by fibre multiplications by unit-length constants. (cf. section 4) Since
$m_\rho$ is $U(1)$-equivariant ($U(1)$ action also defined on $\hat L),\
\widetilde H$  is $U(1)$-invariant, so there exists a $U(1)$-equivariant
diffeomorphism $\phi$ with $\phi_*\hat H = \widetilde H$.  Note that two
choices of such
$\phi$ are different by $U(1)$-equivariant  contact diffeomorphisms, i.e. the
inverse of the one composed with the other belongs to $C^{U(1),0}_{\hat H}$.
Using $\phi$ to pull back the  $U(1)$-invariant CR structure $(\widetilde H,
m_\rho^*\widetilde J|\widetilde H$) on $\hat S$, we obtain a $U(1)$-invariant
CR structure $J=(m_\rho\circ\phi)^*(\widetilde J)|\hat H$ in $\widetilde{\frak
S}^{U(1)}$.  Define $\widetilde\tau: \widetilde P_{ic_0}\to\widetilde{\frak
S}^{U(1)}$ by $\widetilde \tau(\hat L,\widetilde J)=(\hat S,\hat H,J)$ where
$J=(m_\rho\circ\phi)^*(\widetilde J)|\hat H$. The map $\widetilde\tau$ gives
rise to a map
$$
   \tau^t: P_{ic_0}^t\to \frak S^{t,U(1)}.
$$

\n
("uniqueness" of $\|\ \ \|_{\widetilde J}$ by Proposition 2.1) Recall that
$\frak S^{t,U(1)} = \widetilde{\frak S}^{U(1)}/C^{U(1),0}_{\hat H}$. Endow
$\widetilde{\frak S}^{U(1)}$ with the $C^\infty$-topology so that
$\frak S^{t,U(1)}$ has the induced quotient topology.

\bigskip
\n
Proposition 5.2. The map $\tau^t: P_{ic_0}^t\to {\frak S}^{t,U(1)}$
is a homeomorphism onto its image.

\bigskip
\n
Proof. To prove $\tau^t$ is continuous, we will suitably choose a unique
$\|\ \ \|_{\widetilde J}$ and a unique $\phi$ for a given $\widetilde J$.
Remember $\|\ \ \|_{\widetilde J}$ is determined by $\lambda$ in (2.3). We
normalize the solution $\lambda$ of (2.3) by requiring $\lambda=1$ at some
point p, so $\lambda$ is uniquely determined. Furthermore, the map :$\widetilde
J\to\lambda$ is continuous by the standard arguments in the
elliptic theory. (apply the Harnack estimates (Theorems 9.20, 9.22 in [GT])
to get upper bounds for $\log\lambda$ and $\log\lambda^{-1}$ (cf.(3.22))) Let
$Diff_0^{U(1)}(M)$ denote the identity component of $Diff^{U(1)}(M)$. Let
$\frak B_0^{U(1)}$ denote the connected component of $\frak B^{U(1)}$,
containing $\hat H$. To pick up a unique $\phi$, we invoke the following
result.

\bigskip
\n
Lemma 5.3. There is a local smooth tame map $s:\frak B_0^{U(1)}\to
Diff_0^{U(1)}(M)$ near a reference point $H_0$ such that $s(\widetilde
H)_*(H_0)=\widetilde H$.

\bigskip
   We will prove Lemma 5.3 in Appendix A. By Lemma 5.3, the map:$\widetilde
J\to \widetilde H $  composed with $s$ gives a continuous map: $\widetilde
J\to\phi$ near a reference point. We have shown that $\tau^t$ is continuous. On
the other hand, given $J$ in $\widetilde{\frak S}^{U(1)}$, we can extend $J$ to
$\widetilde J$ in $\widetilde P_{ic}$ as below. For $y$ not in the 0-section of
$\hat L$, let $\rho=\|y\|_{\hat J}$ and define $\widetilde J_y$ by the fibre
dilation: $\widetilde J_y(v)= J_x(m_{\rho*}^{-1}(v))$ for $v$ in
$m_{\rho*}\hat H_x, x=m_\rho^{-1}(y)$. Since $J$ is $U(1)$-invariant, it can
be pushed down to define a complex structure $c$ on $\hat N$:
$c(\hat\pi_*(v))=\hat\pi_*(Jv)$ for $v$ in $\hat H_\ell, \ell$ in $\hat S$.
Here we identify $\hat H_\ell$ with the tangent space of $\hat N$
at $\hat\pi(\ell)$. Let $s_0$ denote the 0-section: $\hat N\to\hat L$. For $y$
in $s_0(\hat N)$, we define $\widetilde J_y(v)=s_{0*} c(\hat\pi_*(v))$ for $v$
in $T_y(s_0(\hat N))$. For $v$ tangent to fibres, we just define $\widetilde J$
to be the usual complex structure on $\Bbb C$ in local trivializations. Now it
is a matter to verify that $\widetilde J$ is smooth and hence belongs to
$\widetilde P_{ic}$. First observe that the 2-plane distribution $\frak D$ on
$\hat L$ defined by $m_{\rho*}\hat H (\rho \in \Bbb C\backslash\{0\}$) and
tangent  spaces of $s_0\hat N$ is smooth.   (in a local trivialization $(z,w)$,
this distribution can be described by the kernel of the one-form $ih_zwdz
+ ihdw$. Here $h(z,\bar z)=\|s(z)\|^2_{\hat J}$ for a
local holomorphic section s. cf. (2.4)) To show $\widetilde J$ is smooth, it is
enough to prove $\widetilde J(v)$ is smooth for every smooth vector field $v$.
Write $v=v_{\frak D}+ v_f$. Here $v_{\frak D}$ in $\frak D$ is smooth while
$v_f$ is a smooth vector field tangent to the fibres. It is obvious that
$\widetilde J(v_f)$ is smooth. Let $i_l$ denote the linear isomorphism:
$T_{\hat\pi(\ell)}\hat N\to{\frak D}_\ell$ for $\ell \in \hat L$
so that $i_\ell \circ \hat\pi_*=$ identity on $\frak D_\ell$. Note that
$i_\ell=s_{0*}$ at $\hat\pi(\ell)$ for $\ell$ in $s_0(\hat N)$. Now we can
express
$$
\widetilde J_\ell(v_{\frak D}(\ell))=i_\ell \circ c(\hat\pi_*(v_{\frak
D}(\ell))).
$$

\n
Since $i: \ell\to i_\ell$ is smooth as viewed as a section of
End($\hat\pi^*(T\hat N),\frak D$) over $\hat L$, it follows that $\widetilde
J(v_{\frak D}$) is smooth, hence $\widetilde J$ is smooth. It is not hard
to see that the map $\widetilde{ext}: \widetilde{\frak S}^{U(1)}\to\widetilde
P_{ic}$ defined by $\widetilde{ext}(J)=\widetilde J$ is continuous. (see also
[Ep] for precise construction in a local trivialization and in terms of type
$(0,1)$ vector fields) Moreover, $\widetilde{ext}$  induces a continuous map
$ext$ from ${\frak S}^{t,U(1)}$ to $P_{ic}^t$ by the proof of  Proposition
2.4 and $ext\circ\tau^t$ equals the identity for the same reason. Therefore
$\tau^t$ is injective and hence a homeomorphism onto its image.
            \qquad\qed

\bigskip
   In fact ${\tau}^t$ is surjective onto the connected component of
${\frak S}^{t,U(1)}$, containing the reference element $[\hat J]$. We
will see this below.  First let us determine the universal cover of ($\hat
S,\hat H,\hat J$). Denote the unit disc in the complex plane $C$ by $D$. Define
the hermitian metric $\|\ \ \|_e$ on the trivial holomorphic line bundle $D
\times C$ by
$$
\|(z,w)\|_e = |w|^2/(1-|z|^2)^e
$$

\n
for $e=m/(g-1)$. (recall that $-m$ is the first Chern number of $\hat L$ and
$g$ is the genus of $\hat N$)\ \ It is a direct verification that
$h(z,\bar z)=\|(z,1)\|_e$ satisfies (2.1) in $D$, the universal cover of $\hat
N$. Write an element $A$ in $U(1,1)\times U(1)$ as below:
$$
A = \pmatrix  a&b&0\\    c & d & 0 \\ 0 & 0 & u \endpmatrix
$$

\n
for $u$ in $U(1), \pmatrix a& b\\ c&d\endpmatrix$ in $U(1,1)$ with respect to
the quadratic form given by $\pmatrix 1& 0 \\ 0 & -1\endpmatrix$. The group
$U(1,1)\times U(1)$ acts on $D \times C$ by
$$
A(z,w)=((az+b)/(cz+d), uw/(cz+d)^e).
$$

\n
It is easy to see that $A$ leaves $\|\ \ \|_e$ invariant. (just note that
$|z|^2-|w|^2=|az+bw|^2 - |cz+dw|^2$)  Define $S_e \subset D \times C$ by $\|\ \
\|_e=1.$  It follows that $S_{1/(g-1)}$ is an $m$ to 1 cover of
$S_{m/(g-1)}=S_e$ and a $g-1$ to 1 cover of $S_1.$  Since $S_1\ (=
S^3\backslash\{w=0\})$ is obviously spherical, $S_e$  is spherical too. The
holomorphic line bundle $\hat L$ over the Riemann
surface $\hat N$ gives rise to a representation of $\pi_1(\hat N)$ in
$PU(1,1)\times U(1)$ acting on $D \times C$. Here $PU(1,1)=SU(1,1)/$center
acts on $D$ as holomorphic transformations. It follows that
$S_e/\pi_1(\hat N)$ is a spherical circle bundle of $\hat L$ over $\hat N$ with
the hermitian metric induced from $\|\ \ \|_e$ satisfying (2.1). By uniqueness
(up to a constant multiple) $\hat S = S_e/\pi_1(\hat N)$. As a consequence, the
universal cover of $\hat S$ (as CR manifold), denoted $\widetilde S$,  is the
same for any $(g,c_1)$ and is the infinite cyclic cover of $S_1= S^3\backslash
\{w=0\}$. It is well known (e.g. [BS]) that $\widetilde S$ is homogeneous.

\bigskip
\n
Proposition 5.4. Every element in  $\widetilde{\frak S}^0$ is $U(1)$-invariant
up to a contact diffeomorphism in $C^0_{\hat H}$.

\bigskip
\n
Proof. First note that every spherical CR manifold is locally homogeneous
in the weak sense (i.e. any two points have isomorphic neighborhoods). By
Theorem 8.2 of [ENS] or Theorem 4 of [Go], the universal cover of any element
in $\widetilde{\frak S}^0$ is homogeneous and hence is CR equivalent to
$\widetilde S$ by the classification. ([BS] or [ENS]) Denote $\Gamma$ the
fundamental group of $\hat S.$  It is not hard to see (e.g. [FG] p.44) by the
theorem of Seifert and Van Kampen that $\Gamma$ has a presentation:
$$
  \Gamma = \langle a_1,b_1,...,a_g,b_g,h:\Pi^{i=g}_{i=1}[a_i,b_i]= h^{-m}, h
\ \ \text{central}\rangle
$$

\n
($-m$ is the first Chern number or Euler number). Realize $\Gamma$ as Deck
transformations of $\widetilde S$ via the homomorphism $j:\Gamma\to$
Aut$_{CR}(\widetilde S)$. It is known ([BS], p.234) that
Aut$_{CR}(\widetilde S$) satisfies the following exact sequence:
$$
0\to\Bbb R \to\text{Aut}_{CR}(\widetilde S) \CD @> pr >> \endCD PU(1,1)\to 1.
$$

\n
We claim that $pr \circ j(h) = I$, the identity. Since the quotient space
$j(\Gamma)\backslash\widetilde S$ is CR equivalent to ($\hat S, \hat H, J$) for
some $J$ in $\widetilde{\frak S}^0,$ it is compact, and hence has finite
invariant measure. \ \ Aut$_{CR}(\widetilde S)$ acts on $\widetilde S$
transitively with the compact isotropy group, isomorphic to $U(1)$, so
$j(\Gamma)\backslash$ Aut$_{CR} (\widetilde S)$ is compact, and hence has
finite invariant measure. Let $H = pr
\circ j(\Gamma)$.  It follows that $H\backslash PU(1,1)$ is compact and has
finite invariant measure. By Lemma 5.4 in [Ra], H has property (S) in
$PU(1,1).$  Therefore by Corollary 5.18 in [Ra] the centralizer $Z(H)$ of $H$
in $PU(1,1)$ is the centre of $PU(1,1)$, which consists of the identity. Now
note that $pr \circ j(h)$ is in $Z(H)$ since h is central. Hence $pr \circ j(h)
= I$, so $h$ is mapped into the $\Bbb R$ part by $j$.  Let $dev$ denote the
developing map from $\widetilde S$ onto $S_1 \subset S^3$. Let $hol$ denote
the holonomy map from Aut$_{CR}(\widetilde S)$ onto
Aut$_{CR}(S_1)=PU(1,1)\times U(1) \subset \text{Aut}_{CR}(S^3)= PU(2,1)$. The
developing pair $(hol,dev)$ induces naturally another pair $(hol',
dev'): (\text{Aut}_{CR}(\widetilde S),\widetilde S)\to(\text{Aut}_{CR}
(S_e),S_e)$ by noting that both $S_1$ and $S_e$ are covered by the common
covering space $S_{1/(g-1)}$. Let $a_i', b_i', h'$ denote
the corresponding generators of $a_i, b_i, h$ under the map $hol'\circ j$,
respectively. By projecting the commutator relation in $\Gamma$ into the $U(1)$
part of Aut$_{CR}(S_e)$, we obtain $I=(h')^{-m}$. But for $\hat J, h'=I.$
Hence for $J$ in $\widetilde{\frak S}^0, h'$ is also equal to the identity by
continuity. (note that the representation map $j$ depends on our spherical CR
structure $J$ on $(\hat S, \hat H))$  Thus the subgroup $hol'\circ
j(\Gamma)$ of Aut$_{CR}(S_e)$ can be viewed as a representation of $\pi_1(\hat
N)$ generated by $a_i',b_i'$ in Aut$_{CR}(S_e)$. Therefore
$hol'\circ j(\Gamma)\backslash S_e$ is the spherical circle bundle $S_L$ of
some holomorphic line bundle $L$, determined by Proposition 2.1, as discussed
previously.  (in particular it is $U(1)$-invariant) On the other hand,
$hol'\circ j(\Gamma)\backslash S_e$ is CR equivalent to
$j(\Gamma)\backslash\widetilde S$
representing $(\hat S,\hat H,J)$ in view of $dev'$ being a covering map.
Let $\Sigma$ denote the CR isomorphism from $(\hat S,\hat H,J)$
onto $S_L$. Composing a bundle isomorphism between $\hat L$ and $L$ with a
fibre multiplication map $m_\rho$, we can construct a $U(1)$-equivariant
diffeomorphism $\phi: \hat S\to S_L$.  (note that $\Sigma$ may not be
$U(1)$-equivariant) Now $\phi^* H_L$ and $\phi^*J_L$ are invariant with
respect to the $U(1)$-action on $\hat S.$  By Theorem 5.1 we can find a
$U(1)$-equivariant diffeomorphism $\psi$ such that $\psi\circ
\phi^{-1}\circ \Sigma$ is in $C^0_{\hat H}$ while
$(\psi\circ\phi^{-1}\circ\Sigma)^{-1*}(J)=\psi^{-1*}\circ \phi^*(J_L)$ is
$U(1)$-invariant.  \qquad\qed

\bigskip
   Let $\widetilde{\frak S}^0$ denote the connected component of
$\widetilde{\frak S},$ containing $(\hat S,\hat H,\hat J)$. Let
$\widetilde{\frak S}^{0,U(1)}$ denote the space of $U(1)$-invariant
elements in $\widetilde{\frak S}^0$. Any $(\hat S,\hat H,J)$ in
$\widetilde{\frak S}^{0,U(1)}$ extends to a complex structure $\widetilde
J$ on $\hat L$. The argument in the above proof of Proposition 5.4
shows that $(\hat S,\hat H,J)$ is CR-equivalent to $S_L$ for a certain
holomorphic line bundle $L$. The CR isomorphism between $(\hat S,\hat
H,J)$ and $S_L$ implies the existence of a holomorphic
bundle isomorphism between $(\hat L, \widetilde J)$ and $L$ in view of the
proof of Proposition 2.4. Since $S_L$ is uniquely determined by $L$
(Proposition 2.1), it follows that $(\hat S,\hat H,J)$ is uniquely
determined by $(\hat L,\widetilde J)$, i.e. suppose two $(\hat S,\hat
H,J_i)$ have isomorphic extensions ($\hat L, \widetilde J_i),\ i=1,2,$ then
($\hat S,\hat H,J_1$) is CR-equivalent to ($\hat S,\hat H,J_2$).
Furthermore by Proposition 3.14 in [Ep] we have

\bigskip
\n
 Lemma 5.5. Let ($\hat L,\widetilde J_i$) be the extension
of ($\hat S,\hat H,J_i)$ in $\widetilde{\frak S}^{0,U(1)}, i=1,2$. Suppose
($\hat L,\widetilde J_1$) is isomorphic to ($\hat L,\widetilde J_2$) by a
bundle  automorphism  in $Bdiff_0$. Then ($\hat S,\hat H,J_1)$ is CR-equivalent
to ($\hat S,\hat H,J_2)$ by a contact diffeomorphism in $C^{U(1),0}_{\hat H}$.

\bigskip
   Let $\widetilde{\frak S}^{U(1),0}$ denote the connected component of
$\widetilde{\frak S}^{U(1)},$ containing $(\hat S,\hat H,\hat J$).
Define ${\frak S}^t_0$ to be the quotient space of
$\widetilde{\frak S}^{U(1),0}$ modulo $C^{U(1),0}_{\hat H}$ (or $C^0_{\hat H}$:
two quotient spaces are the same by the above discussion), i.e. two elements in
$\widetilde{\frak S}^{U(1),0}$ are equivalent if one is carried to another
by a contact diffeomorphism in $C^{U(1),0}_{\hat H}$ (or $C^0_{\hat H}$
resp.) by pulling back.  Observe that ${\frak S}^t_0$ is an open connected
subset of $\widetilde{\frak S}^{0,U(1)}/C^0_{\hat H}$ which
equals $\widetilde{\frak S}^0/C^0_{\hat H}$ in view of Proposition 5.4. Since a
CR equivalence $\phi$ between two $U(1)$-invariant CR circle bundles
is $U(1)$-equivariant,  we have ${\frak S}^{t,U(1)}\ (=\widetilde
{\frak S}^{U(1)}/C^{U(1),0}_{\hat H})=\widetilde{\frak S}^{U(1)}/C^0_{\hat
H}$.

\bigskip
\n
Proposition 5.6. $\tau^t: P_{ic_0}^t\to{\frak S}^t_0$ is
surjective and a homeomorphism in view of Proposition 5.2.

\bigskip
\n
Proof. Given an element ($\hat S,\hat H,J$) in $\widetilde{\frak S}^{U(1),0}$,
there associates an extension ($\hat L,\widetilde J$) in
$\widetilde P_{ic_0}$. We claim $\tau^t([(\hat L,\widetilde J)])=[(\hat
S,\hat H,J)].$  Recall that the construction of $\tau^t$ involves a map
$m_\rho$ and a $U(1)$-equivariant diffeomorphism $\phi$ on $\hat S$.
Extend $\phi$ to a bundle automorphism $\widetilde\phi$ in $Bdiff_0.\ \
\widetilde\tau (\hat L,\widetilde J)=(\hat S,\hat H,(m_\rho\circ\phi)^*(\hat
J)|\hat H)$ is the restriction of $(m_\rho\circ\widetilde\phi)^*(\widetilde J)$
on $\hat L$ to $(\hat S,\hat H).$  Since $m_\rho\circ\widetilde\phi$ is a
bundle automorphism of $\hat L$ in $Bdiff_0$, it follows that ($\hat
S,\hat H,(m_\rho\circ\phi)^*(\widetilde J))$ is CR-equivalent to ($\hat
S,\hat H,J)$ by a contact diffeomorphism in $C^{U(1),0}_{\hat H}$ according
to Lemma 5.5.          \qquad\qed

\bigskip
   We remark that in [KT] Kamishima and Tan studied the deformation space of
$U(1)$-invariant spherical CR-structures by analyzing the space of developing
pairs. Their deformation space for $M = \hat S$ is in one-to-one
correspondence with our space ${\frak S}^{t,U(1)}$ by "contact"
reduction. According to Corollary 5.2.2 in [KT], this space is homeomorphic to
Hom($\pi_1(\hat N),PU(1,1))/PU(1,1) \times T^{2g}$, and it is well known
that the dimension of Hom($\pi_1(\hat N),PU(1,1))/PU(1,1)$ is $6g-6$, the
dimension of Teichmuller space.(e.g. [Go]) Thus the total dimension is
$6g-6 + 2g = 8g-6$ (cf. Theorem 6 (d) in [Go]) while Proposition 5.6 shows that
an open connected subset ${\frak S}^t_0$ of ${\frak S}^{t,U(1)}$ is
homeomorphic to $P_{ic_0}^t$ of the same dimension by Theorem 3.9.

   Next we want to endow $\frak S^t_0$ with a natural differentiable
structure through the general local slice theorem, and with this differentiable
structure on $\frak S^t_0, \tau^t$ in Proposition 5.6 is a diffeomorphism.
Given $J$ in $\widetilde{\frak S}^{U(1),0}$, there passes a local slice $\frak
S$ of $\frak J$ according to Theorem A of [CL2]. Let $\frak P$ denote the
diffeomorphism given in Theorem A (1) of [CL2]. Define $\psi: \widetilde
P_{ic_0}\to\frak S$ near ($\hat L,\widetilde J$) with $\widetilde \tau\ (\hat
L,\widetilde J) = (\hat S,\hat H,J)$ by
$$
 \psi = \text{proj}_S \circ \frak P^{-1}\circ\widetilde\tau.   \tag 5.1
$$

\n
Here proj$_S$ denotes the projection onto the $\frak S$-component. Since the
pullback by a contact diffeomorphism does not change the vanishing of the
Cartan tensor, $\psi$ actually maps into $\frak Q^{-1}(0)$. At $(\hat
L,\widetilde J),$ there passes a local slice, denoted $\frak S_{Pic}$, by
Lemma 3.7. We claim $\psi: \frak S_{Pic}\to\frak Q^{-1}(0)$ is an immersion
(between
two finite dimensional manifolds) by choosing unique $\rho$ and $\phi$ in
defining $\tau$ as explained in the proof of Proposition 5.2. First note that
the action of bundle automorphisms does not change the transversality of
tangent vectors at $\widetilde J$ in $\widetilde P_{ic_0}$. (here
transversality means transverse to the orbit of $\frak B$ or $Bdiff_0$ acting
on $\widetilde J$) Use the bundle automorphism $m_\rho\circ\widetilde\phi$ in
the proof of Proposition 5.6 to reduce our immersion problem to the following:

\bigskip
\n
Lemma 5.7. Let $J$ in $\widetilde {\frak S}^{U(1),0}$ be the restriction of its
extension $\widetilde J$ in  $\widetilde P_{ic_0}$. Let $\widetilde J'$ be an
infinitesimal variation of $\widetilde J$ and $J'$ be the
corresponding infinitesimal variation of $J$ in $\frak J$. Suppose $J'$ is
tangent to the orbit of $C^0_{\hat H}$ acting on $J$. Then $\widetilde J'$ is
also tangent to the orbit of $\frak B$ acting on $\widetilde J.$

\bigskip
\n
Proof. Take a local trivialization $(z,w)$ of $\hat L$ relative to $\widetilde
J$ with $w$ being the fibre coordinate. Let $s$ be the local holomorphic
section given by $z\to (z,1).$  Let $h = h(z,\bar z)=\|s(z)\|^2$ where
the hermitian metric $\|\ \ \|$ is chosen according to Proposition 2.1. By
Lemma  5.5 $\hat S$ is precisely discribed by $\|\ \ \|=1$ or $hw\bar w =1$ in
the above local trivialization. It is easy to verify that $Z=\pa_z - (\log h)_z
w\pa_w$ is tangent to $\hat S$. Let $\theta^1= dz, \theta^2=dw +(\log h)_z
wdz$. Then $\{\theta^1,\theta^2\}$ is dual to $\{Z,\pa_w\}$. Now we
can write $J = i\theta^1\otimes Z +$ conjugate and $\widetilde J= J +
(i\theta^2\otimes\pa_w +$conjugate). Moreover let $\widetilde J_t$ be a family
of extensions of $J_t$ with $\widetilde J_0=\widetilde J, J_0=J$. Let $Z_t= Z +
a_t \bar Z$ be a frame of type (1,0) with respect to $J_t$ with
$a_0=0$. (i.e. an eigenvector of $J_t$ with eigenvalue $i$) Let $\theta^1_t
 = (\theta^1 - (\bar a_t)\theta^{\bar 1})/(1-|a_t|^2)$. It is straightforward
to determine $\theta^2_t$ such that $\{\theta^1_t,\theta^2_t,
\theta^{\bar 1}_t,\theta^{\bar 2}_t\}$ is dual to $\{Z_t, \pa_w, \bar
Z_t, \pa_{\bar w}\}$. The result is $\theta^2_t = dw + (\log h)_z
w\theta^1_t+\bar a_t (\log h)_z w\theta^{\bar 1}_t.$

It follows that $\widetilde J_t = J_t + (i\theta^2_t \otimes\pa_w$ +
conjugate). Computing the derivative at $t=0$ gives
$$
   \widetilde J_t'= J_t'= 2ia_t'dz \otimes\bar Z +\ \text{conjugate}
$$

\n
by observing that $\theta^{2'}_t=0$. Writing $\widetilde J_t'= E_1{}^{\bar
1}dz\otimes\pa_{\bar z} + E_1{}^* \bar w dz \otimes\pa_{\bar w}$ + conjugate
(cf. (3.6)),  we obtain  $E_1{}^{\bar 1}= 2ia_t', E_1{}^*=-2ia_t'(\log
h)_{\bar z}.$  Hence $E_1=E_1{}^* + E_1{}^{\bar 1}\bar\Gamma =0$ (cf.(3.24)) by
noting that  $\Gamma=(\log h)_z$. Now by the assumption we can write $J'=
B'_J(f)= f_{,1}{}^{\bar 1}\theta^1\otimes\bar Z$ + conjugate.  (in [CL1]
we write $Z_1$ instead of $Z$, and choosing the specific contact form (2.4),
we have the torsion $A_1{}^{\bar 1}$ to vanish by (2.6)) That is to say,
$2ia_t'= f_{,1}{}^{\bar 1}$. To show $\widetilde J'= P(V)$ for some $V$
represented by $(v^1,v)$ (cf. (3.31)), we take $v^1= f,{}^1/2i$ and $v=-\mu
f/2i$. It follows that $P(V)=2i(v^1{}_{,\bar 1}, v_{,\bar 1} + \mu v_{\bar
1})=(f_,{}^1{}_{\bar 1}, 0)=(E_{\bar 1}{}^1, 0$) which represents $\widetilde
J'.$     \qquad\qed

\bigskip
\n
Proof of Theorem A: By Lemma 5.7 the differential of $\widetilde\tau$  maps a
tangent vector of $\frak S_{Pic}$ at $\widetilde J$ to a tangent vector
transverse to the orbit of $C^0_{\hat H}$ acting on $J$. It follows that the
differential of $\psi$ (cf.(5.1)): $\frak S_{Pic}\to\frak Q^{-1}(0)$ is
injective.  Therefore $\psi|\frak S_{Pic}$ is an immersion by the inverse
function theorem, so $\psi|\frak S_{Pic}$ gives rise to a local coordinate map
for $\frak S^t_0$ near [J] in view of Proposition  5.6. Transition functions of
these coordinate maps are smooth because $\frak S_{Pic}$'s can be viewed as
local coordinate neighborhoods for the smooth manifold $P_{ic}^t$. (cf. Theorem
C) Thus $\frak S^t_0$ is a smooth manifold and in the way to define its
differentiable structure we actually obtain that $\tau^t$ is a diffeomorphism.
Hence the dimension of $\frak S^t_0$ equals the dimension of $P_{ic}^t$, which
is $2(4g-3)=8g-6$ by Theorem C.   \qquad\qed

\bigskip
\n
Proof of Corollary $B: {\iota}$ is well defined in view of Proposition 5.4.
>From the proof of Proposition 5.4 we learn that any spherical $(\hat
S,\hat H,J)$ in $\widetilde{\frak S}^0$ is CR-equivalent to $S_L$ for some
holomorphic line bundle L. Let $\phi: (\hat S,\hat H,J)\to(S_L,H_L,J_L)$
denote this CR-isomorphism. Take $\theta=\phi^*(\theta_L)$. It is
obvious that $(\hat S,\hat H,J,\theta)$ is in $\widetilde{\frak
M}^0_{-1,0}$ and $\iota$ maps [$(\hat S,\hat H,J,\theta$)] to
[($\hat S,\hat H,J$)]. Thus $\iota$ is surjective and hence bijective in view
of Corollary 2.3.
On the other hand it is easy to see that both $\iota$ and its inverse are
continuous, so $\iota$ is a homeomorphism.     \qquad\qed

\bigskip
\n
Appendix A: The $U(1)$-invariant version of Gray's theorem

   We will prove Theorem 5.1 and Lemma 5.3. First denote the smooth $U(1)$
action by $U_\rho, 0 \le \rho \le 2\pi$ with $U_0 = U_{2\pi}$. Pick a
$U(1)$-invariant metric $g$. (which can be obtained by averaging the action on
an arbitrary metric) Any $U(1)$-invariant contact bundle $H$ in $\frak
B^{U(1)}$ can be uniquely determined by a $U(1)$-invariant 1-form $\theta$ with
$\theta\land d\theta\ne 0$ and $|\theta|_g=1$. Here $|\ |_g$ denotes the
pointwise length with respect to the metric $g$. Still denote the space of all
such 1-forms by $\frak B^{U(1)}.$

\bigskip
\n
Lemma A.1. $\frak B^{U(1)}$ is a tame Frechet submanifold of the tame Frechet
manifold $\frak B$.

\bigskip
\n
Proof: Let $\Omega^1 \ (\Omega^1_{U(1)}$, respectively) denote the space of all
smooth ($U(1)$-invariant, respectively) 1-forms on our closed manifold M. It is
known that $\Omega^1$ is a tame Frechet space ([Ha]). Since the process of
averaging the $U(1)$ action on a 1-form is a tame linear map from $\Omega^1$ to
$\Omega^1_{U(1)}$, it follows that $\Omega^1_{U(1)}$ is a tame direct summand,
hence a tame Frechet space. (Lemma 1.3.3 on p.136 in [Ha]) Now consider the
space
$$
  T_\theta\frak B^{U(1)}:= \{\text{smooth 1-form }\ \eta:\langle
           \eta,\theta\rangle_g=0, U_\rho^*\eta=\eta\}
$$

\n
where $\langle , \rangle_g$ denote the pointwise inner product with respect to
$g$. It is easy to see that the linear map $proj: \Omega^1_{U(1)}\to
T_\theta\frak B^{U(1)}$ given by $proj(\eta)=\eta-\langle\eta,\theta\rangle_g
\theta$ is tame. Therefore $T_\theta\frak B^{U(1)}$ is a tame direct summand of
$\Omega^1_{U(1)}$, hence a tame Frechet space.  Define a map $\Phi_\theta:
T_\theta\frak B^{U(1)}\to\frak B^{U(1)}$ by
$$
  \Phi_\theta(\eta) = (\eta + \theta)/|\eta + \theta|_g.
$$

\n
If we endow $\frak B^{U(1)}$ with the $C^\infty$ topology, then $\Phi_\theta$
is a local homeomorphism near 0 with its inverse $\pi_\theta$ given by
$$
   \pi_\theta(\eta) = \eta/\langle \eta,\theta\rangle_g -\theta.
$$

\n
Now we can compute the transition function for the overlap of two
neighborhoods centered at $\theta$ and $\theta'$:
$$
  \pi_{\theta'}\circ \Phi_\theta(\eta)=(\eta+\theta) /
\langle\eta+\theta,\theta'\rangle_g -\theta'.
$$

\n
It is easy to see that $\pi_{\theta'}\circ \Phi_\theta$ is smooth tame. We
have shown that $\frak B^{U(1)}$ is a tame Frechet manifold. Actually the map
$\Phi_\theta$ also parametrizes $\frak B$ near $\theta$. Therefore $\frak
B^{U(1)}$ is a tame Frechet submanifold of $\frak B$.     \qquad\qed

\bigskip
\n
Lemma A.2. $Diff^{U(1)}(M)$ is a tame Frechet submanifold of $Diff(M)$, the
group of smooth diffeomorphisms on $M$. Moreover, $Diff^{U(1)}(M)$ is a smooth
tame Lie group.

\bigskip
\n
Proof: Let $T_{e} Diff(M)(T_e Diff^{U(1)}(M)$, respectively) denote
the space of all smooth ($U(1)$-invariant, respectively) vector fields on $M$.
Here $e$ denotes the identity diffeomorphism. It is easy to see that the map
$pr: T_e Diff(M)\to T_e Diff^{U(1)}(M)$ given by
$$
   pr(X) = \f 1{2\pi}\int_0^{2\pi}U_{\rho*}(X)d\rho
$$

\n
is linear and tame. Therefore $T_e Diff^{U(1)}(M)$ is a tame direct
summand of the tame Frechet space $T_e Diff(M)$. It follows that
$T_e Diff^{U(1)}(M)$ is also a tame Frechet space. Given a smooth vector
field $X$ on $M$, we denote $\exp_pX$ the time =1 point of the geodesic (with
respect to the invariant metric $g$) passing through $p$ in $M$ with the
velocity $X$. Define $\Psi: T_e Diff(M)\to Diff(M)$ by $\Psi(X)(p)=
\exp_{p}X.$  It is known that $\Psi$ parametrizes $Diff(M)$ near $e$.
Moreover, $\Psi$ maps the subspace $T_e Diff^{U(1)}(M)$ (injectively for sure)
into $Diff^{U(1)}(M)$ since $U_\rho$'s are isometries with respect to $g$. We
claim that $\Psi$ restricted to $T_e Diff^{U(1)}(M)$ is actually surjective
onto $Diff^{U(1)}(M)$. Suppose $\phi$ is a $U(1)$-equivariant diffeomorphism
near $e$ and let $X= \Psi^{-1}(\phi)$. We need to show that $X$ is
$U(1)$-invariant.  Let $\phi_t$ denote the geodesic flow of $X$ with respect to
$g$. Since $U_\rho$ is an isometry (with respect to $g$),\ $U_\rho \circ
\phi_t(p)$ is a geodesic connecting $U_\rho(p)\ (t=0)$ and $U_\rho \circ
\phi(p)(t=1)$ for $p$ in $M$. On the other hand, $\phi_t \circ U_\rho(p)$ is a
geodesic connecting $U_\rho(p)\ (t=0)$ and $\phi \circ U_\rho(p)(t=1)$ which
equals $U_\rho \circ \phi(p)$ by the assumption. Now for $\phi$ close enough to
$e$, the uniqueness of geodesics connecting two points in a convex neighborhood
implies $U_\rho \circ \phi_t = \phi_t \circ U_\rho$. It follows that $U_{\rho*}
\circ X = X \circ U_\rho,$ i.e. $X$ is $U(1)$-invariant, so $\Psi$ parametrizes
$Diff^{U(1)}(M)$ as a tame Frechet submanifold of $Diff(M)$ near $e$. For a
$U(1)$-equivariant diffeomorphism $\psi\ne e$, the local parametrization
$\Psi_\psi$ defined by $\Psi_\psi(X)= \Psi(X)\circ\psi$ for $Diff(M)$ also
parametrizes $Diff^{U(1)}(M)$ near $\psi$ when $X$'s are restricted to
$T_e Diff^{U(1)}(M)$. We have shown that $Diff^{U(1)}(M)$ is
a tame Frechet submanifold of $Diff(M)$. It is easy to see that
$Diff^{U(1)}(M)$ is a group under composition and $Diff(M)$ is a smooth tame
Lie  group (p.148 in [Ha]). It follows that $Diff^{U(1)}(M)$ is also a
smooth tame Lie group.  \qquad\qed

\bigskip
\n
PROOF OF THEOREM 5.1: The action $P: Diff^{U(1)}(M)\times \frak B^{U(1)}
\to \frak B^{U(1)}$ is described by
$$
   P(\phi, \theta)=\phi^*\theta/|\phi^*\theta|_g.
$$

\n
It is easy to see that the maps:($\phi, \theta)\to\phi^*\theta$ and
$\eta\to \langle\eta,\eta\rangle_g^{\f 12}=|\eta|_g$ are smooth tame, so $P$ is
smooth tame. Let $b_\theta(X)= D_1P(e,\theta)(X)$, the partial derivative with
respect to the first variable of $P$ at the identity $e$. Let $H$ denote the
contact bundle annihilated by the contact form $\theta$. Let $\phi_t$ be a
smooth family of $U(1)$-equivariant diffeomorphisms such that $\phi_0 = e$ and
$\f d{dt}|_{t=0}\phi_t = X$. We compute
$$
\align
       \f d{dt}|_{t=0}|\phi_t^*\theta|_g^{-1}
    &= \f d{dt}|_{t=0}(|\phi_t^*\theta|_g^2)^{\f{-1}2}\\
    &= (-\f 12)\cdot 2\langle X\rfloor  d\theta, \theta\rangle_g\\
    &= -\langle X\rfloor  d\theta, \theta\rangle_g\endalign
$$

\n
for $X$ tangent to $H$. It follows that for $X$ tangent to $H$,
$$
\align
   b_\theta(X)& = \f d{dt}|_{t=0} P(\phi_t,\theta)\\
              & = \f d{dt}|_{t=0} \phi_t^*\theta +
                  \f d{dt}|_{t=0}|\phi_t^*\theta|_g^{-1}\theta\\
              & = X \rfloor d\theta - \langle X  \rfloor d\theta,
                    \theta\rangle_g\theta\\
              & = \pi_{H^*}(X \rfloor d\theta).\endalign
$$

\n
Here $\pi_{H^*}$ is the projection onto the orthogonal complement of $\theta$.
Now given $\eta$ in $T_\theta\frak B^{U(1)}$, we want to find $X$ such that
$\pi_{H^*}(X \rfloor d\theta)=\eta$. It is easy to find a unique $X$ tangent to
$H$ so that $X\rfloor d\theta = \eta$ on $H$.  Since $\theta, H$, and
$\eta$ are all $U(1)$-invariant, it follows that $X$ is also $U(1)$-invariant
by uniqueness. On the other hand, $X\rfloor d\theta=\pi_{H^*}(X \rfloor
d\theta)$ on $H$. But $\eta$ is orthogonal to $\theta$. Thus
$b_\theta(X)=\pi_{H^*}(X \rfloor d\theta)=\eta.$  We have proved that the map:
$X\to b_\theta(X)=D_1P(e,\theta)(X)$ from $T_e Diff^{U(1)}(M)$ to
$T_\theta\frak B^{U(1)}$ is surjective with a right inverse $\eta\to X$. It is
easy to check that the linear map:$\eta \to X$ is tame. Now our theorem follows
from Theorem 2.4.1 on p.198 in [Ha].    \qquad\qed

\bigskip
\n
PROOF OF LEMMA 5.3: First observe that in the proof of Theorem 2.4.1 on
p.198 in [Ha], we actually show that the action with a reference point fixed
is locally surjective and has a smooth tame right inverse by Theorem 1.1.3 on
p.172 in [Ha]. This means in our case the action
$\phi\to P(\phi,\theta_0)$ with $\theta_0$ fixed is locally surjective and
has a smooth tame right inverse $V$. Set $s(\widetilde
H)=(V(\widetilde\theta))^{-1}$ where $\widetilde\theta$ is the contact form
associated to the contact bundle $\widetilde H$ near $H_0$.
    \qquad\qed

\bigskip
\n
Appendix B: \ \ An infinitesimal slice of $\widetilde{\frak M}_{-1,0}/C_{\hat
H}$

Take  a family of pseudohermitian structures $(J_{(t)}, \theta_{(t)})$ on
$(\hat S, \hat H)$ with $J_{(0)} = \hat
J, \theta_{(0)} = \hat \theta$.  At $t = 0$, express
$$
(\overset \,\cdot \to J_{(t)},\overset \,\cdot \to \theta_{(t)})=(2E_1{}^{\bar
1} \hat\theta^1\otimes \hat Z_{\bar 1} + 2E_{\bar 1}{}^1 \hat\theta^{\bar 1}
\otimes \hat Z_1, 2h\hat \theta)    \tag B.1
$$

\n
where $E_1{}^{\bar 1}$ is a deformation tensor at $\hat J$ and $h$ is just a
real-valued function.  (see (2.14) on p.231 in [CL1]; also note
$\theta_{(t)}|\hat H = 0$) \ \ Next we observe the action of $C_{\hat
H}$. Let $\phi_t \in C_{\hat H}$ be a family of contact diffeomorphisms
with $\phi_0 =$ identity.  Compute
$$
\aligned
\f d{dt}|_{t=0}(\phi^*_t\hat J, \phi^*_t\hat\theta) &= (L_{X_f}\hat J,
L_{X_f}\hat\theta) \ \ \text{(Lemma 3.4 on p.239 in [CL1])}\\
 &= (2B'_{\hat J}f, -(\hat Tf)\hat \theta) \ \ \text{((3.13) and the proof of
Lemma 3.4 in [CL1])} \endaligned
$$

\n
where $B'_j$ is the second-order operator defined on p.236 in [CL1] and $\hat
T$ is the vector field uniquely determined by $\hat \theta(\hat T) = 1, \hat T
\lrcorner d\hat \theta = 0.$  Define
$$
\widetilde B'_{\hat J}f = B'_{\hat J}f - \f 12(\hat T f)\hat \theta.
$$

\n
Then we have the following orthogonal decomposition:
$$
T_{(\hat J, \hat \theta)}\{(J, \theta): \theta|\hat H = 0\} = Ker\widetilde
B_{\hat J}\oplus\ \text{Range}\ \widetilde B'_{\hat J}     \tag B.2
$$

\n
where $\widetilde B_{\hat J}$ is the adjoint operator of $\widetilde B'_{\hat
J}$, given by $$
\widetilde B_{\hat J}(\widetilde E) = B_{\hat J}E + \f 12 h_{,o}
$$

\n
for $\widetilde E = E + h\hat\theta, E = E_1{}^{\bar 1}\hat\theta^1
\otimes \hat Z_{\bar 1} + E_{\bar 1}{}^1 \hat\theta^{\bar
1}\otimes \hat Z_1$.  Here $B_{\hat J}$ is defined on p.235 in [CL1].  Note
that Range $\widetilde B'_{\hat J}$ is the tangent space of the orbit of $
C_{\hat H}$ passing through $(\hat J, \hat \theta)$.  The decomposition (B.2)
is valid either in $L^2$ category or in $C^\infty$ category mainly because of
the fourth-order operator $\widetilde B_{\hat J}\widetilde B'_{\hat J} =
\Delta_{\hat J} - \f
14\hat T^2 = \f 12\frak L^*_\al\frak L_\al + O_2 \ (\al = i\sqrt{\f 72})$ being
subelliptic.  (see Lemma 2.1 in [CL2])  Now linearizing the equations
$R_{J,\theta} \equiv -1, \ A_{J,\theta}\equiv 0$ at $(\hat J, \hat\theta)$ in
the direction $(\overset\,\cdot\to J, \overset\,\cdot\to\theta)$ given by
(B.1), we obtain
$$
\cases
& i(E_1{}^{\bar 1}{}_{,\bar 1}{}^1 - E_{\bar 1}{}^1{}_{,1}{}^{\bar 1}) + 2h +
  4\Delta_bh = 0\\
& E_{\bar 1}{}^1{}_{,0} + 2h_{,\bar 1}{}^1 = 0   \endcases \tag B.3
$$

\n
by (5.15), (5.9) in [Lee] and (2.20), (2.18) in [CL1], where $\Delta_bh =
-(h_{,1}{}^1 + h_{,\bar 1}{}^{\bar 1})$.  (see (4.10) in [Lee])  Since
elements in Range $\widetilde B'_{\hat J}$ satisfy the linear equations (B.3),
we get from (B.2) that  an infinitesimal slice of $\widetilde{\frak M}
_{-1,0}/C_{\hat H}$  is the intersection of $Ker \widetilde B_{\hat J}$ and the
solution space of (B.3). Write $\widetilde K = K + k\theta$ in this
infinitesimal slice.   It follows  that $\widetilde K$ satisfies the following
system of equations:
$$
\cases
& K_1{}^{\bar 1}{}_{,\bar 1}{}^1 + K_{\bar 1}{}^1{}_{,1}{}^{\bar 1} + \f
12k_{,0} = 0\\
& i(K_1{}^{\bar 1}{}_{,\bar 1}{}^1 - K_{\bar 1}{}^1{}_{,1}{}^{\bar 1}) + 2k +
4\Delta_bk = 0\\
& K_{\bar 1}{}^1{}_{,0} + 2k_{,\bar 1}{}^1 = 0.   \endcases
$$

\medskip
\midinsert
\vspace{.5in}
\endinsert
\noindent
$$
\left. \aligned
&\text{Cheng: Institute of Mathematics}\\
&      \text{Academia Sinica, Nankang}\\
&      \text{Taipei, Taiwan, R.O.C.}\\
 &\text{E-mail: Cheng\@math.sinica.edu.tw}
\endaligned
\qquad\qquad
\aligned
&\text{Tsai: Department of Mathematics}\\
 &      \text{National Taiwan University}\\
  &     \text{Taipei, Taiwan, R.O.C.}\\
&\text{E-mail: ihtsai\@math.ntu.edu.tw}
\endaligned
\right.
$$
\vfill\eject\end